\documentclass[11pt]{article}
\def\mytitle{Conic relaxation approaches for equal deployment problems}

\usepackage{amssymb}
\usepackage{amsthm}
\usepackage{color}
\usepackage{enumerate}
\usepackage{url}
\usepackage{graphicx}
\usepackage{amsmath}  
\usepackage{multirow}
\usepackage{lineno, blindtext}
\usepackage{changepage}
\usepackage[nottoc]{tocbibind}

\def\@themcountersep{}
\newtheorem{THEO}{Theorem}[section]
\newtheorem{ALGO}[THEO]{Algorithm}
\newtheorem{ASSUM}[THEO]{Assumption}

\newtheorem{LEMM}[THEO]{Lemma}

\newtheorem{REMA}[THEO]{Remark}

\def\0{\mbox{\bf 0}}
\def\1{\mbox{\bf 1}}
\def\2{\mbox{\bf 2}}
\def\3{\mbox{\bf 3}}
\def\4{\mbox{\bf 4}}
\def\5{\mbox{\bf 5}}
\def\6{\mbox{\bf 6}}
\def\7{\mbox{\bf 7}}
\def\8{\mbox{\bf 8}}
\def\9{\mbox{\bf 9}}

\def\b{\mbox{\boldmath $b$}}
\def\c{\mbox{\boldmath $c$}}

\newdimen\zhige \zhige=0pt
\def\chige#1{{\setbox\zhige\hbox{#1}\ifdim\ht\zhige=1ex\accent24 #1%
  \else\ooalign{\unhbox\zhige\crcr\hidewidth\char24\hidewidth}\fi}}

\def\e{\mbox{\boldmath $e$}}

\def\g{\mbox{\boldmath $g$}}
\def\h{\mbox{\boldmath $h$}}

\def\l{\mbox{\boldmath $l$}}

\def\u{\mbox{\boldmath $u$}}
\def\v{\mbox{\boldmath $v$}}

\def\x{\mbox{\boldmath $x$}}
\def\y{\mbox{\boldmath $y$}}
\def\z{\mbox{\boldmath $z$}}
\def\A{\mbox{\boldmath $A$}}
\def\B{\mbox{\boldmath $B$}}
\def\C{\mbox{\boldmath $C$}}

\def\F{\mbox{\boldmath $F$}}

\def\H{\mbox{\boldmath $H$}}

\def\O{\mbox{\boldmath $O$}}

\def\V{\mbox{\boldmath $V$}}

\def\X{\mbox{\boldmath $X$}}
\def\Y{\mbox{\boldmath $Y$}}

\def\FC{\mbox{$\cal F$}}

\def\IC{\mbox{$\cal I$}}

\def\KC{\mbox{$\cal K$}}

\def\NC{\mbox{$\cal N$}}

\def\SC{\mbox{$\cal S$}}

\def\Real{\mbox{$\mathbb{R}$}}

\def\SMAT{\mbox{$\mathbb{S}$}}

\ifx11 % 11 => use these commands, 10 => NOT use these commands

\fi

\newcommand{\futoji}[1]{\mbox{\boldmath $#1$}}

\begin{document}

\noindent \textcolor{black}{\leaders\hrule width 0pt height 0.08cm \hfill}
\vspace{0.5cm} 

%{\Large \bf \noindent Fast computation for optimal selection problem in tree breeding 
%using Second-Order Cone Programs}
{\Large \bf \noindent \mytitle}

\vspace{0.2cm}
 \noindent
Sena Safarina
\footnote[1]{
Department of Mathematical and Computing Science,
 Tokyo Institute of Technology, 2-12-1-W8-29 Ookayama, Meguro-ku, Tokyo
 152-8552, Japan.
 (The work of M. Yamashita was partially supported by JSPS KAKENHI (Grant-in-Aid for Scientific Research (C), 15K00032).)
 },
%Pietro Belotti
%\footnote[2]{
%Xpress Optimizer team, FICO, Starley Way, Birmingham, B37 7GN, UK. 
%},
%Mituhiro Fukuda
%\footnotemark[1],
Satoko Moriguchi
\footnote[2]{
% BioSylve, PO Box 30709, Nairobi 00100, Kenya, \newline
School of Business Administration,
Faculty of Urban Liberal Arts,
Tokyo Metropolitan University,
1-1 Minami-Osawa, Hachioji-shi, Tokyo 192-0397, Japan.
},
Tim J. Mullin
\footnote[3]{
% BioSylve, PO Box 30709, Nairobi 00100, Kenya, \newline
The Swedish Forestry Research Institute (Skogforsk), Box 3, S\"{a}var 918 21, Sweden; 
and 224 rue du Grand-Royal Est, QC, J2M 1R5, Canada.
(The work of T. J. Mullin was partially supported by F\"{o}reningen Skogstr\"{a}dsf\"{o}r\"dling (The Swedish Tree Breeding Foundation)},
and
Makoto Yamashita
\footnotemark[1]
% (Corresponding Author) 
% Department of Mathematical and Computing Science,
% Tokyo Institute of Technology, 2-12-1-W8-29 Ookayama, Meguro-ku, Tokyo
% 152-8552, Japan (Makoto.Yamashita@is.titech.ac.jp).
\\
Submitted: February 27, 2017.
% Submitted: March 1, 2015.

% \noindent \textcolor{blue}{\leaders\hrule width 0pt height 0.02cm \hfill}
\vspace{0.3cm}

\noindent {\bf Abstract:} 
% {\bf Oh! This is a very good paper.}
%\begin{linenumbers}
An important problem in the breeding of livestock, crops, and forest trees is the optimum of selection of genotypes that maximizes genetic gain.
The key constraint in the optimal selection is a convex quadratic constraint that 
ensures genetic diversity, therefore, the optimal selection can be cast as a second-order cone programming (SOCP) problem.
Yamashita et al. (2015) exploits the structural sparsity of the quadratic constraints
and reduces the computation time drastically while attaining the same optimal solution.

This paper is concerned with the special case of equal deployment (ED), in which we solve the optimal selection problem with the constraint that contribution of genotypes must either be a fixed size or zero.
This involves a nature of combinatorial optimization, and the ED problem can be described as a
mixed-integer SOCP problem.

In this paper, we discuss conic relaxation approaches for the ED problem based on LP (linear programming), 
SOCP, and SDP (semidefinite programming). We analyze theoretical bounds derived
from the SDP relaxation approaches using the work of Tseng (2003)
and show that the theoretical bounds are not quite sharp for tree breeding problems.
We propose a steepest-ascent method that
combines the solution obtained from the conic relaxation problems
with a concept from discrete convex optimization
in order to acquire 
an approximate solution for the ED problem in a practical time.
From numerical tests, we observed that among the LP, SOCP, and SDP relaxation problems,
SOCP gave a suitable solution from the viewpoints of the optimal values 
and the computation time.
The steepest-ascent method starting from the SOCP solution provides high-quality solutions much faster than 
an existing method that has been widely used for the optimal selection problems
and a branch-and-bound method.
%\end{linenumbers}
\vspace{0.3cm} 


\noindent {\bf MCS2010 classification:} 
90C05  	Linear programming,
90C11  	Mixed integer programming,
90C22  	Semidefinite programming,
90C25  	Convex programming,
90C59  	Approximation methods and heuristics,
90C90  	Applications of mathematical programming,
92-08  Biology and other natural sciences (Computational methods).

\noindent \textcolor{black}{\leaders\hrule width 0pt height 0.08cm \hfill}
\vspace{0.5cm} 

\section{Introduction}

Computational methods based on mathematical optimization have started gaining attention from breeding researchers,
since the optimization methods provide efficient approaches and give theoretical aspects for 
the optimality of the obtained solutions.
For example, optimal selection problems that determine the contributions of genotypes
are studied for clonal seed orchards
and dairy cattle~\cite[etc]{bomosssoul1993optimal, hallander2009optimum, hinrichs2006algorithm, lindgren1989deployment, MEUWISSEN97, mullin2014opsel}.
% dairy cattle is from hinrichs2006algorithm % This includes very large-scale data.

A main objective in optimal selection problems is to attain the highest response from a genotype selection.
Lindgren et al.~\cite{lindgren1989deployment} proposed a linear deployment in which the genotype contributions
are basically proportional to their breeding values. This deployment was derived from a concept that the genotypes with 
higher breeding values should appear more frequently than those with lower values.
An advantage of the linear deployment was the extremely low computation cost, 
since it could be computed by a greedy algorithm.
However, the linear deployment worked well only when the pedigree situation was simple,
that is, the candidate genotypes were unrelated.
If the selected genotypes do not embrace enough diversity, 
the response will critically diminish through inbreeding depression~\cite{charlesworth1987inbreeding, williams1996inbreeding}
due to accumulated kinship.

Meuwissen~\cite{MEUWISSEN97} introduced a quadratic constraint to control 
a group coancestry under an appropriate level.
%Mathematical optimization problems have become
%important mathematical tool for breeding researchers 
%to determine an optimal contributions of candidate pedigree members.
%Therefore, it is natural that 
%computation methods based on mathematical optimization started 
%gaining attention from breeding researchers, in particular,  the optimization methods 
%provide efficient approaches and give clues to 
%the optimality of the obtained  solution.
He developed the Lagrangian multiplier method to maximize the genetic response 
with the quadratic constraints.
This method was implemented in a software package GENCONT~\cite{MEUWISSEN97}, 
and it has been widely accepted among breeding researchers. A serious drawback of the Lagrangian multiplier method is that 
this method does not always generate optimal solutions.
In contrast, Pong-Wong et al.~\cite{PONG-WONG07} employed an SDP approach. 
This approach is based on mathematical optimization, and 
% when the contributions are continuous variables.
they demonstrated that this approach gave the optimal contributions exactly.
This approach was extended in~\cite{ahlinder2014using}, but their SDP approach required long computation time
even when they used parallel computing with the help of 
SDPA (a high-performance solver for SDPs)~\cite{SDP-HANDBOOK, SDPA6}.
Recently, Yamashita et al.~\cite{yamashita2015efficient} proposed 
an SOCP (second-order cone programming) approach
and successfully reduced the computation time of the SDP approach
attaining the same optimal solution.

The problems solved by the SDP approach~\cite{PONG-WONG07} and 
the SOCP approach~\cite{yamashita2015efficient} 
are unequal deployment (UD) problems of form
% A mathematical form of the UD problem is given as follow.
\begin{eqnarray}
% (UD) &:& 
\begin{array}{lcl}
\max &:& \g^T \x \\
\mbox{subject to} &:& \x^T \A \x \le 2 \theta, \\
& & \e^T \x = 1, \\
& & \l \le \x \le \u.
\end{array} \label{eq:UD}
\end{eqnarray}
Throughout this paper, we use $Z$ to denote the number of candidate genotypes. In the UD problem, the variable is the vector $\x \in \Real^Z$, and $x_i$ indicates the contribution
of the $i$th genotype.
We use  a superscript $T$ to denote the transpose of a vector or a matrix.
The cost vector $\g \in \Real^Z$ in the objective function is 
the estimated breeding value (EBV)~\cite{LYNCH98}.
Since this vector is computed separately, we regard $\g$ as a constant vector.
The matrix $\A \in \Real^{Z \times Z}$ is the Wright numerator matrix~\cite{WRIGHT22}.
The elements of this matrix are given from the information of heredity diagram.
We should emphasize that the matrix $\A$ is always symmetric and positive definite.
Hence, with a given constant $\theta > 0$, the constraint
$\x^T \A \x \le 2 \theta$ is a convex constraint, 
and this quadratic constraint ensures that the group coancestry $\frac{\x^T \A \x}{2}$ in the selected group
is kept under a permissible range $\theta$.
We use  $\e \in \Real^n$ to denote the vector of all ones,
therefore, the constraint $\e^T \x = 1$ indicates that the total contribution of
all the candidates is unity.
In addition, the vectors $\l \in \Real^Z$ and $\u \in \Real^Z$ are the lower and upper bounds
of the variable $\x$, respectively.

The name an \textit{unequal} deployment indicates that 
the contributions need not to be equal.
Since the variable $\x$ is a continuous variable and the constraints are linear or convex-quadratic,
 the UD problem can be cast a type of SOCP problems, as pointed in \cite{yamashita2015efficient}.
Therefore, the UD problem can be solved in a polynomial time algorithm,
for example, interior-point methods for SOCP~\cite{ALIZADEH03, DOMAHIDI13, TSUCHIYA99}.

This paper is concerned with the special-case problem of \textit{equal} deployment (ED) form
\begin{eqnarray}
\begin{array}{lclcl}
OPT_{ED} &:=& \max &:& \g^T \x \\
& & \mbox{subject to} &:& \x^T \A \x \le 2 \theta, \\
& & & & \e^T \x = 1, \\
& & & & \l \le \x \le \u, \\
& & & & x_1, \ldots, x_n \in \left\{0, \frac{1}{N}\right\}.
\end{array} \label{eq:ED}
\end{eqnarray}
We use $OPT_{ED}$ to denote the optimal value of this problem.
The crucial difference from the UD problem is that the ED problem has the binary constraints
$x_1, \ldots, x_n \in \left\{0, \frac{1}{N}\right\}$.
We choose exactly $N$ genotypes from $Z$ candidates, and the 
selected $N$ genotypes must contribute their genes equally.
The ED problems fit breeding populations, where we consider
the selected genotypes should contribute with the same amount
and therefore we 
require a fixed-size population. 

Weng et al.~\cite{weng2012unequal} solved the ED problem 
only with the linear constraints and the binary constraints using the ``Solver'' tool in Microsoft Excel. Meuwisen extended GENCONT to the ED problems incorporating some heuristic methods
so that GENCONT generated approximate solutions that satisfy the binary constraints.
The heuristic methods implemented in GENCONT
are partially discussed in \cite{woolliams2015genetic}.

% the detail of the heuristic method implemented in GENCONT
% is not unapparent, since its source code is not open.

From the viewpoint of mathematical optimization, the most difficult constraint 
 $\x^T \A \x \le 2 \theta$  is a quadratic convex constraint. An ED problem (\ref{eq:ED}) can thus be viewed as a mixed-integer second-order cone
programing (MI-SOCP) problem.
%From the binary constraint, the number of feasible points 
%we search for an optimal solution is bounded.
%Therefore, the appearance of the ED problem is much easier than the UD problem.
%On the contrary to this intuition, this binary constraint makes the ED problem much 
%harder than the UD problem. 
%Actually, MI-SOCP problems are considered as NP-hard problems.
%
Many approaches have been explored to solve MI-SOCP efficiently.
Ben-tal and Nemirovski~\cite{ben2001polyhedral} proposed a polyhedral relaxation that approximates
a second-order cone with a polyhedron so that the resulting problem can be handled
with software packages for mixed-integer linear programming.
Drewes applied an outer approximation method and a branch-and-cut method~\cite{drewes2012subgradient}.
For other approaches, a survey paper due to
Benson and Saglam~\cite{benson2013mixed} is a good reference.

Theoretically speaking, 
MI-SOCP is an SOCP problem with integer constraints, hence,
we can obtain an exact optimal solution if we rely on 
the branch-and-bound framework.
However, we suffer from a long computation time if we pursue the exact solution.
For example, CPLEX can directly handle MI-SOCP problems, but fails to complete the computation of a small case $Z = 1050$ and $N = 50$
(it tried to choose $N = 50$ genotypes from $Z = 1050$ candidates) in one week.
% due to a structure of the group coancestry constraint.
Mullin and Belotti~\cite{mullin2016using} combined the outer approximation 
method and the branch-and-bound method and reduced the computation time.
However, it also requires half a day for the small case $Z = 200$ to attain 
the gap $0.5\%$, so it is still hard to say 
that this approach is practical for larger instances $Z \ge 5000$.
To manage ED problems in a practical time, it is desirable that we find 
a high-quality approximate solution instead of the exact solution. 

In this paper, we propose an integration of conic relaxation approaches
and a steep-ascent method originally developed for discrete convex functions
to derive a suitable solution for practical usage in a reasonable computation time. % for the ED problems.

An epoch-making paper on conic relaxation approach was the application of SDP problem to the max-cut problems
by Goemans and Williamson~\cite{GOEMANS95}. They converted a feasible set of 
the max-cut problems into the space of positive semidefinite matrices with the rank-one constraint on the matrix variable,
and they derived an SDP problem by ignoring this rank-one constraint. They showed that 
a solution generated with a randomized algorithm from an optimal solution of the resulting SDP problem
gave very good approximation to the original max-cut problem.
%The feasible set $\FC_{\mbox{mc}}$ of the max-cut problems 
%can be expressed with binary constraints; $\FC_{\mbox{mc}} := \{\x \in \Real^n: x_i = \pm 1 \mbox{ for } i=1,\ldots, n\}$.
%They embedded $\FC_{\mbox{mc}}$ into the space of symmetric matrices.
%$\FC_{\mbox{MC}} := \{\X \in \SMAT^n : X_{ii} = 1 \mbox{ for } i=1,\ldots, n, \ \X \succeq \O, \ \mbox{rank}(\X) = 1\}$.
%Throughout this paper, we use $\SMAT^n$ for the space of $n \times n$ symmetric matrices and 
%the notation $\X \succeq \O$ indicates a symmetric $\X$ is positive semidefinite.
Following this achievement, the SDP relaxation approach has widely been applied to combinatorial optimization problems,
see \cite{wolkowicz2012handbook} and the references therein.
Theoretical evaluation of the quality of the approximate solution were discussed in
\cite[etc]{he2008semidefinite, hsia2015improved, nesterov1998semidefinite, tseng2003further, ye1999approximating}.
Conic relaxation approaches are the relaxation approaches that employs linear programming (LP), 
SOCP or SDP problems. A remarkable points of the three conic programming problems (LP, SOCP, and SDP)
 is that they can be analyzed in the framework of Euclidean Jordan 
algebras~\cite{faybusovich1997euclidean, schmieta2001associative}. Hence, 
the resulting relaxation problems can be solved in polynomial time by interior-point methods~\cite{nesterov1994interior}
and many software packages are available~\cite{sturm1999using, TODD99, SDP-HANDBOOK}.
Kim and Kojima~\cite{kim2001second} reported a numerical evaluation on the relaxation approaches using LP, SOCP, and SDP  for some quadratic optimization problems.

%On the other hand, discrete convex optimization has another abundant research direction.
%We might consider that a convex function in continuous space is a discrete convex function
%if we restrict the variable space to the integer points.
%However, this intuition is not always true.
%In discrete convex 
%optimization~\cite{murota2003discrete},
%the concept of convexity is divided into two convex concepts, L-convex and M-convex.
%If a function is an M-convex function,
%a steep-descent method proposed in~\cite{murota2004steepest} can find 
%its global minimizer.

On the other hand, discrete convex optimization has another abundant research direction. We might consider that a convex function in continuous space is a discrete convex function if we restrict the variable space to the integer points, although this naive intuition is not appropriate because such a function does not always have useful properties of convex functions, and some deep combinatorial or discrete-mathematical considerations are needed for discrete convexity. In the theory of discrete convex analysis~\cite{murota2003discrete}, two convexity concepts, called L-convexity and M-convexity, play primary roles. L-convex functions and M-convex functions are convex functions with additional combinatorial properties distinguished by "L" and "M", which are conjugate to each other through a discrete version of the Legendre-Fenchel transformation. If a function is an M-convex function, a step-descent method proposed in~\cite{murota2004steepest} can find its global minimum.

In this paper, we first introduce conic relaxation problems for the ED problems,
and discuss the relations between the relaxation problems.
We analyze the theoretical bounds of the randomized algorithm starting from
the solution of the SDP relaxation problem.
However, when we numerically evaluate these bounds using tree-breeding datasets,
we learn that these bounds are not so sharp.
Instead of pursuing an exact solution by branch-and-bound frameworks
that impose heavy computation costs, 
our focus is to acquire a favorable solution that is available in a practical
computation time.
To obtain such a solution, 
we develop a steep-ascent method that employs the solution obtained from
the conic relaxation problems as a starting point.
The usual steep-descent method~\cite{murota2004steepest} minimizes an objective function 
on a particular feasible set. Since the ED problem
is a maximization problem, we consider a steep-ascent method instead of 
a steep-descent method.
We embed the quadratic constraint $\x^T \A \x \le 2 \theta$
into the objective function as a penalty term with a weight computed from
the Lagrange multiplier.
This new objective function is not an M-concave function,
therefore, we cannot guarantee that the solution obtained by 
the steep-ascent method is a global solution of the ED problem.
However, through numerical experiments, 
we observe that the steep-ascent method generates qualified solutions
for the ED problem.
In particular, the steep-ascent method starting with the SOCP relaxation problem
attains the best performance among the LP, SOCP, and SDP relaxation problems.
% the trade-off between the solution quality and the computation time
%is a reasonable approach for the ED problems.
%Among the conic relaxation problems, the SDP relaxation gives the tightest
%approximation. The quality of the SOCP
%relaxation is competitive to that of the SDP relaxation, but the computation time
%of the SOCP relaxation is much less than that of the SDP relaxation.
%When we consider the trade-off between the solution quality and the computation time,
%the steep-descent method starting with the SOCP relaxation problem
%is a reasonable approach for the ED problems.
Actually, we verify from numerical experiments that 
this approach performs better than existing methods like GENCONT in the viewpoints of both
solution quality and computation time.

The rest of this paper is organized as follows.
In Section~2, we introduce LP, SOCP, and SDP relaxation problems for the ED problems,
and we discuss the strength of these conic relaxations.
In Section~3, we analyze the approximation rate of the SDP relaxation 
based on the work of Tseng~\cite{tseng2003further}.
Section~4 gives the details of the steep-ascent method
specialized for the ED problems.
In Section~5, we present numerical results to compare the conic relaxations
and to evaluate the solution acquired by the steep-ascent method.
We also compare this result with existing methods.
In Section~6, we will give a conclusion and discuss future directions.

\subsection{Notation}
We use $|S|$ to denote the cardinality of a set $S$.
The vector $\e_S$ is the vector of all ones of the lengths $|S|$.
In contrast, we denote by $\e_i$ the vector of all zeros except one in the $i$th position.
The symbol $\SMAT^n$ is used to denote the space of $n \times n$ symmetric matrices,
and $\X \succeq \O$ indicates that a symmetric matrix $\X$ is positive semidefinite.
The inner-product between $\A \in \SMAT^n$ and $\X \in \SMAT^n$
is defined by $\A \bullet \X := \sum_{i=1}^n \sum_{j=1}^n A_{ij} X_{ij}$.
The trace of  a matrix $\A \in \SMAT^n$ is given by $\mbox{Trace}(\A) := \sum_{i=1}^n A_{ii}$. For a vector $\x \in \Real^n$, $\x \ge \0$ indicates the element-wise non-negativity of 
$\x$, that is, $x_1, \ldots, x_n \ge 0$.

\section{Conic relaxations for equally deployment problems}

In this section, we first derive an SDP relaxation problem of an ED problem. 
Then, by a further relaxation
of the positive semidefinite condition using a relaxation
technique proposed in \cite{kim2001second},
we obtain an LP relaxation problem.
Finally, we apply a continuous relaxation technique to the ED problem to obtain 
an SOCP relaxation problem. The reason we employ a different relaxation approach
for only the SOCP relaxation is that we can exploit a structural sparsity 
in the Wright numerator matrix $\A$.

A standard form of SDP problems can be given as follow:
\begin{eqnarray}
% (SDP) : 
\begin{array}{lcl}
\min &:& \C \bullet \X \\
\mbox{subject to} &:& \F_i \bullet \X =  b_i \ (i=1,\ldots,m),\\
% \mbox{subject to} &:& \F_i \bullet \X (\le, =, \ge)  b_i \ (i=1,\ldots,m)\\
& & \X \succeq \O.
\end{array}\label{eq:SDP}
\end{eqnarray}

In this standard form, the variable matrix is $\X \in S^n$. 
% The symbol $(\le, = \ge)$ can 
% take either $\le$, $=$, or $\ge$ for each $i$.
% that  $\F_i \bullet \X (\le, =, \ge)  b_i$ takes one of 
% $\F_i \bullet \X \le  b_i$, $\F_i \bullet \X =  b_i$ or $\F_i \bullet \X \ge b_i$ for each $i$.
The input matrices in (\ref{eq:SDP}) are $\C, $ $\F_1, \ldots,\F_m \in \SMAT^n$,
while the vector $\b \in \Real^n$ is an input vector.
Shortly speaking, a standard SDP form 
% in (\ref{eq:SDP}), we
minimizes a linear objective function
over linear constraints and a positive semidefinite condition on $\X$.

As a first step to derive an SDP relaxation from the ED problem (\ref{eq:ED}),
we remove the variables that can be fixed from the box constraints.
More precisely, if $l_i > 0$, we fix $x_i = \frac{1}{N}$. Similarly, we fix $x_i = 0$
if $u_i < \frac{1}{N}$. We ignore the cases $l_i > \frac{1}{N}$,
$u_i < 0$ or $l_i > u_i$,
since we can immediately detect the infeasibility of the ED problem.
Then, we define two sets $F$ and $V$
so that the two sets separate the set $\{1, \ldots, Z\}$ disjointly
and $x_i$ is fixed to $c_i \in \left\{0,\frac{1}{N}\right\}$ for $i \in F$ 
while $x_i$ remains as a decision variable for $i \in V$.

Without loss of generality, we assume that $V = \{1, 2, \ldots,  |V|\}$,
$F = \{ |V| + 1, |V| + 2, \ldots, Z\}$, and $g_1 \ge g_2 \ge \ldots \ge g_{|V|}$.
Along with these $V$ and $F$, we introduce the vectors $\x_V$ and $\c_F$
that divide $\x \in \Real^Z$ into the two parts
$\x = \left(\begin{array}{cc} \x_V \\ \c_F \end{array}\right)$.
We also divide the Wright numerator matrix $\A$ into the four parts;
$\A = \left(\begin{array}{cc} \A_{VV} & \A_{VF} \\ 
\A_{FV} & \A_{FF} \end{array}\right)$.
The sizes of $\A_{VV}$, $\A_{FV}( = \A_{VF}^T)$, and $\A_{FF}$ are $|V| \times |V|$,
$|F| \times |V|$, and $|F| \times |F|$, respectively.
We further partition the vectors and the matrices that appear in the ED problem
into the corresponding parts;
\begin{eqnarray}
\begin{array}{lclcl}
OPT_{ED} & = & \max &:& \g_V^T \x_V + \g_F^T \c_F\\
& & \mbox{subject to} &:& \x_V^T \A_{VV} \x_V + 2 \c_F^T \A_{FV} \x_V
+ \c_F^T \A_{FF} \c_F 
\le 2 \theta ,\\
& & & & \e_V^T \x_V +  \e_F^T \c_F = 1,  \\
& & & & x_i \in \left\{0, \frac{1}{N}\right\} \mbox{ for } i \in V.
\end{array} \label{eq:ED2}
\end{eqnarray}
Note that we also removed the box constraints $\l \le \x \le \u$ from the ED problem
by fixing the variables in $\x_F$ to $\c_F$. 
We count the number of $x_i$ that is fixed to $c_i$
by $p:= \left| \left\{ i \in F : x_i = \frac{1}{N} \right\} \right|$.
Therefore, we will choose $N-p$ genotypes from $|V|$ candidates  in (\ref{eq:ED2}), while 
we choose $N$ genotypes from $Z$ candidates in the original ED problem~(\ref{eq:ED}).
\begin{REMA}\label{re:remarkPN}
We can assume $p \le N$ and $|V| \ge 2$ without loss of generality.
In the case $p > N$, we can detect the infeasibility of the problem~(\ref{eq:ED}).
If $|V| = 1$, we have $F = \{2, \ldots, Z\}$. Therefore, $x_1$ is also fixed with $x_1 = 1-\sum_{i=2}^Z c_i$,
and all the variables can be fixed without solving~(\ref{eq:ED2}).
\end{REMA}

We change
the decision variables by $\y_V := 2 N \x_V - \e_V \in \Real^{|V|}$
and we use $y_i$ to denote the $i$th element of $\y_V$.
% $x_1, \ldots, x_{|V|}$
% into $y_1, \ldots, y_{|V|}$ with the relation $y_i = 2 N x_i - 1$ 
Then, the binary constraints 
$x_1, \ldots, x_{|V|} \in \left\{0, \frac{1}{N}\right\}$
are mapped to $y_1, \ldots, y_{|V|} \in \left\{-1, 1 \right\}$.
Even without employing this variable change, 
we can also directly apply the SDP relaxation method in a similar way to \cite{gorge2012semidefinite}.
The reason we employed this variable change 
is for the later discussion in Section~4 so that 
most of the matrices $\B^k$ there will be diagonal matrices.

We will denote the $i$th element of $\y_V$ by $y_i$.
We define $g_{\min} := \min \{g_i : i = 1, \ldots, Z\}$, 
$\bar{\g}_V: = \frac{1}{4N} (\g_V - g_{\min} \e_V)$, 
$\bar{g} := \frac{1}{2N} (\g_V - g_{\min} \e_V)^T \e_V 
+ (\g_F - g_{\min} \e_F)^T \c_F + g_{\min}$,
 $\bar{\c}_F := \A_{VV} \e_V + 2 N \A_{VF} \c_F$,
 $\bar{\theta} := 2N^2(2\theta - \c_F^T \A_{FF} \c_F) 
- \frac{1}{2}\e_V^T \A_{VV} \e_V - 2N \c_F^T \A_{FV} \e_V$,
and $\bar{N} := 2N(1-\e_F^T \c_F) - |V| = 2(N - p) - |V|$.
From these definitions, it is easy to check $\bar{\g}_V \ge \0$
and $\g^T \x = 2 \bar{\g}_V^T \y_V + \bar{g}$ using $\e^T \x = 1$. 
We now have another expression of the ED problem;
\begin{eqnarray}
\begin{array}{lclcl}
OPT_{ED} &=& \max &:&  2 \bar{\g}_V^T \y_V 
+ \bar{g} \\
& & \mbox{subject to} &:& \y_V^T \A_{VV} \y_V + 2 \bar{\c}_F^T  \y_V 
\le 2 \bar{\theta},\\
& & & & \e_V^T \y_V = \bar{N},\\
& & & & y_i \in \left\{-1, 1\right\} \mbox{ for } i \in V. 
\end{array} \label{eq:ED3}
\end{eqnarray}
% We removed a constant term $yyy$ from the objective function.

% to this expression.
By introducing a variable matrix $\Y_{VV} \in \SMAT^{|V|}$, 
we apply the lift-and-project method of Lov{\'a}sz and 
Schrijver~\cite{lovasz1991cones}. As a result,
we obtain one more equivalent form;
\begin{eqnarray}
\begin{array}{lclcl}
OPT_{ED} &=& \max &:&  
\left(\begin{array}{cc} 0 & \bar{\g}_V^T \\ \bar{\g}_V & \O \end{array}\right)
\bullet
\left(\begin{array}{cc} 1 & \y_V^T \\ \y_V & \Y_{VV} \end{array}\right)
+ \bar{g} \\
& & \mbox{subject to} &:& \left(\begin{array}{cc} -2\bar{\theta} & \bar{\c}_F^T \\ \bar{\c}_F & \A_{VV} \end{array}\right)
\bullet
\left(\begin{array}{cc} 1 & \y_V^T \\ \y_V & \Y_{VV} \end{array}\right)
\le 0,\\
& & & & \left(\begin{array}{cc} -2\bar{N} & \e_V^T \\ \e_V & \O \end{array}\right)
\bullet
\left(\begin{array}{cc} 1 & \y_V^T \\ \y_V & \Y_{VV} \end{array}\right)
= 0,\\
& & & & \left(\begin{array}{cc} -\bar{N^2} & \0^T \\ \0 & \e_V \e_V^T  \end{array}\right)
\bullet
\left(\begin{array}{cc} 1 & \y_V^T \\ \y_V & \Y_{VV} \end{array}\right)
= 0,\\
% & & Y_{ii} = 1 \mbox{ for } i = 1, \ldots, |V| \\
& & & & \left(\begin{array}{cc} -1 & \0^T \\ \0 & \e_i \e_i^T  \end{array}\right)
\bullet
\left(\begin{array}{cc} 1 & \y_V^T \\ \y_V & \Y_{VV} \end{array}\right)
= 0   \mbox{ for } i \in V,\\
& & & & \left(\begin{array}{cc} 1 & \y_V^T \\ \y_V & \Y_{VV} \end{array}\right)
\succeq \O, \quad 
\mbox{rank}
\left( \left(\begin{array}{cc} 1 & \y_V^T \\ \y_V & \Y_{VV} \end{array}\right) \right)
= 1.
\end{array} \label{eq:ED4}
\end{eqnarray}
The key property for the equivalence between (\ref{eq:ED3}) and \ref{eq:ED4}
is $\Y_{VV} = \y_V \y_V^ T$ from the rank-1 constraint on the matrix $
\left(\begin{array}{cc} 1 & \y_V^T \\ \y_V & \Y_{VV} \end{array}\right)$.
We will denote the $(i,j)$th element of $\Y_{VV}$ by $Y_{ij}$.
The equality $Y_{ii} = y_i^2$ for $i = 1, \ldots,  |V|$ should holds
for feasible solution of (\ref{eq:ED4}),
hence
$ \left(\begin{array}{cc} -1 & \0^T \\ \0 & \e_i \e_i^T  \end{array}\right)
\bullet
\left(\begin{array}{cc} 1 & \y_V^T \\ \y_V & \Y_{VV} \end{array}\right)
= 0$ leads to the binary constraint $y_i \in \{-1, 1\}$.
In (\ref{eq:ED4}), 
we introduced a redundant constraint $(\e_V\e_V^T) \bullet \Y_{VV} = \bar{N}^2$
that was derived from $(\e_V^T \y_V)^2 = \bar{N}^2$ and $\Y_{VV} = \y_V \y_V^T$.
It is known that redundant constraints of this type make the SDP relaxation tighter,
and we can often obtain better approximate solution.
The hardest constraint in (\ref{eq:ED4}) is the rank-1 constraint.
This constraint embraces a nature of combinatorial optimization.
By removing this hardest constraint,
we build an SDP relaxation problem and we denote its optimal value by 
$OPT_{SDP}$.

\begin{eqnarray}
\begin{array}{lclcl}
OPT_{SDP} &:=& \max &:&  
\left(\begin{array}{cc} 0 & \bar{\g}_V^T \\ \bar{\g}_V & \O \end{array}\right)
\bullet
\left(\begin{array}{cc} 1 & \y_V^T \\ \y_V & \Y_{VV} \end{array}\right)
+ \bar{g} \\
& & \mbox{subject to} &:& \left(\begin{array}{cc} -2\bar{\theta} & \bar{\c}_F^T \\ \bar{\c}_F & \A_{VV} \end{array}\right)
\bullet
\left(\begin{array}{cc} 1 & \y_V^T \\ \y_V & \Y_{VV} \end{array}\right)
\le 0,\\
& & & & \left(\begin{array}{cc} -2\bar{N} & \e_V^T \\ \e_V & \O \end{array}\right)
\bullet
\left(\begin{array}{cc} 1 & \y_V^T \\ \y_V & \Y_{VV} \end{array}\right)
= 0 ,\\
& & & & \left(\begin{array}{cc} -\bar{N^2} & \0^T \\ \0 & \e_V \e_V^T  \end{array}\right)
\bullet
\left(\begin{array}{cc} 1 & \y_V^T \\ \y_V & \Y_{VV} \end{array}\right)
= 0,\\
% & & Y_{ii} = 1 \mbox{ for } i = 1, \ldots, |V| \\
& & & & \left(\begin{array}{cc} -1 & \0^T \\ \0 & \e_i \e_i^T  \end{array}\right)
\bullet
\left(\begin{array}{cc} 1 & \y_V^T \\ \y_V & \Y_{VV} \end{array}\right)
= 0   \mbox{ for } i \in V,\\
& & & & \left(\begin{array}{cc} 1 & \y_V^T \\ \y_V & \Y_{VV} \end{array}\right)
\succeq \O.
\end{array} \label{eq:SDP-R}
\end{eqnarray}

%Kim and Kojima~\cite{kim2001second} introduced an SOCP relaxation for non-convex
%quadratic optimization problems via SDP relaxation techniques. 
%A matrix $\X \in \SMAT^{n}$ is positive semidefinite if $\u^T \X \u \ge 0$ for $\forall \u \in \Real^n$,
%and a second-order cone of dimension $q$ is defined by 
%$\KC^q = \{ \x \in \Real^q : x_1 \ge \sqrt{\sum_{i = 2}^n x_i^2} \}$.
%They first used the Schur complement to convert 
%the positive semidefinite constraint
%of the lift-and-projection method 
%$\left(\begin{array}{cc} 1 & \y^T \\ \y & \Y \end{array}\right) \succeq \O$ 
%into an equivalent condition $\Y -\y\y^T \succeq \O$.
%Then, they chose a finite number of vectors $\u_1, \ldots, \u_K$ to relax 
%this condition  into second-order cone constraints;
%\begin{eqnarray*}
%& & \Y - \y\y^T \succeq \O 
% \Leftrightarrow 
%\u^T (\Y - \y\y^T) \u \ge 0 \mbox{ for } \forall \u \in \Real^n\\
%& \overset{\mbox{(relaxation)}}\Rightarrow & 
%% \left(\begin{array}{cc} 1 \\ \u_i \end{array}\right)^T \left(\begin{array}{cc} 1 & \y^T \\ \y & \Y \end{array}\right) 
%% \left(\begin{array}{cc} 1 \\ \u_i \end{array}\right) \ge 0 \mbox{ for } i = 1, \ldots, k
%\u_i^T (\Y - \y\y^T) \u_i \ge 0 \mbox{ for } i = 1, \ldots, k 
% \Leftrightarrow 
%\left(\begin{array}{cc} (\u_i \u_i^T)\bullet \Y + 1 \\ (\u_i \u_i^T)\bullet \Y - 1 \\ \u_i^T \y \end{array}\right)
%\in \KC^3 \mbox{ for } i = 1, \ldots, K. \\
%\end{eqnarray*}
When we further relax the positive semidefinite constraint
$\left(\begin{array}{cc} 1 & \y_V^T \\ \y_V & \Y_{VV} \end{array}\right)
\succeq \O$, 
we can obtain an LP relaxation problem.
In general, a matrix $\X \in \SMAT^{n}$ is positive semidefinite if and only if 
$\u^T \X \u \ge 0$ for $\forall \u \in \Real^n$.
For the positive semidefinite constraint of (\ref{eq:SDP-R}),
we choose a set of vectors 
$\u_{ij} = \e_i - \e_j \in \Real^{1+|V|}$ for $i = 1,\ldots, |V|$ and $j = i+1, \ldots, |V|+1$
as a subset of $\Real^{1+|V|}$.
We use $\hat{W}$ to denote 
the non-diagonal upper-triangular position of $\Y_{VV}$, that is 
$\hat{W}:= \{(i,j) \in V \times V :  i < j\}$.
The key step to derive an LP relaxation problem is the following step:
\begin{eqnarray*}
& & \left(\begin{array}{cc} 1 & \y_V^T \\ \y_V & \Y_{VV} \end{array}\right) \succeq \O \\
&\Leftrightarrow & 
\u^T \left(\begin{array}{cc} 1 & \y_V^T \\ \y_V & \Y_{VV} \end{array}\right) \u \ge 0
 \mbox{ for }  \forall \u \in \Real^{1+|V|} \\
&\overset{\mbox{(relaxation)}}\Rightarrow & 
\u_{ij}^T \left(\begin{array}{cc} 1 & \y_V^T \\ \y_V & \Y_{VV} \end{array}\right) \u_{ij} \ge 0
 \mbox{ for }  i = 1, \ldots, |V| \mbox{ and } j = i+1, \ldots, |V|+1 \\
& \Leftrightarrow & 
\left\{\begin{array}{lll}
Y_{ii} \ge y_i^2 & \mbox{for} & i \in V, \\
Y_{ii} Y_{jj} \ge Y_{ij}^2 & \mbox{for} & (i,j) \in \hat{W}.
\end{array}\right.
\end{eqnarray*}
From the constraints
$Y_{ii} = 1$ for $ i \in V$
% $i = 1, \ldots, |V|$ 
in (\ref{eq:SDP-R}), the constraints 
$Y_{ii} \ge y_i^2$ and $Y_{ii} Y_{jj} \ge Y_{ij}^2$ 
are linear constraints in nature.
Consequently, we reach an LP relaxation problem, whose optimal value 
is denoted as $OPT_{LP}$.

\begin{eqnarray}
\begin{array}{lclcl}
OPT_{LP} &=& \max &:&  
\left(\begin{array}{cc} 0 & \bar{\g}_V^T \\ \bar{\g}_V & \O \end{array}\right)
\bullet
\left(\begin{array}{cc} 1 & \y_V^T \\ \y_V & \Y_{VV} \end{array}\right)
+ \bar{g} \\
& & \mbox{subject to} &:& \left(\begin{array}{cc} -2\bar{\theta} & \bar{\c}_F^T \\ \bar{\c}_F & \A_{VV} \end{array}\right)
\bullet
\left(\begin{array}{cc} 1 & \y_V^T \\ \y_V & \Y_{VV} \end{array}\right)
\le 0,\\
& & & & \left(\begin{array}{cc} -2\bar{N} & \e_V^T \\ \e_V & \O \end{array}\right)
\bullet
\left(\begin{array}{cc} 1 & \y_V^T \\ \y_V & \Y_{VV} \end{array}\right)
= 0 ,\\
& & & & \left(\begin{array}{cc} -\bar{N^2} & \0^T \\ \0 & \e_V \e_V^T  \end{array}\right)
\bullet
\left(\begin{array}{cc} 1 & \y_V^T \\ \y_V & \Y_{VV} \end{array}\right)
= 0,\\
% & & Y_{ii} = 1 \mbox{ for } i = 1, \ldots, |V| \\
& & & & \left(\begin{array}{cc} -1 & \0^T \\ \0 & \e_i \e_i^T  \end{array}\right)
\bullet
\left(\begin{array}{cc} 1 & \y_V^T \\ \y_V & \Y_{VV} \end{array}\right)
= 0  \mbox{ for } i \in V,\\
& & & & \begin{array}{lcl} 
-1 \le y_i \le 1  &\mbox{ for }&  i \in V,\\
 -1  \le Y_{ij} \le 1  &\mbox{ for }&  (i,j)  \in \hat{W},
\end{array}\\
& & & & \Y_{VV} \in \SMAT^{|V|}.
\end{array} \label{eq:LP-R}
\end{eqnarray}

We now move our focus to an SOCP relaxation problem. 
In a similar way to the above step that derives (\ref{eq:LP-R}) from (\ref{eq:SDP-R}),
it may be possible to apply an SOCP relaxation technique developed in~\cite{kim2001second} to (\ref{eq:SDP-R}).
In contrast, we utilize a continuous relaxation technique 
that converts the binary constraint $x_i \in \left\{0, \frac{1}{N}\right\}$ into a continuous constraint
% $x_i \in \{ x_i : 0 \le x_i \le \frac{1}{N}\}$.
$0 \le x_i \le \frac{1}{N}$.
The main reason of this continuous relaxation is 
that we can keep the efficient SOCP formula 
of~\cite{yamashita2015efficient} that 
extensively exploits a structural sparsity of the Wright numerator matrix $\A$.

A second-order cone of dimension $q$ is defined by 
$\KC^q :=\left \{ \x \in \Real^q : x_1 \ge \sqrt{\sum_{i = 2}^n x_i^2} \right\}$.
A standard form of second-order cone programming (SOCP) problem in this paper is given as follows:
\begin{eqnarray}
% (SOCP) : 
\begin{array}{lcl}
\max &:& \c^T \x \\
\mbox{subject to} &:& \F  \x =  \b,\\
% \mbox{subject to} &:& \F  \x (\le, =, \ge)  \b \\
& & \h - \H \x \in \KC^{q}.
\end{array}\label{eq:SOCP}
\end{eqnarray}

The decision variable here is $\x \in \Real^n$ and the objective function is a linear function
with a constant vector $\c \in \Real^n$.
The linear constraints are encoded with a matrix $\F \in \Real^{m \times n}$ and a vector $\b \in \Real^m$.
The second-oder cone constraint is given with a vector $\h \in \Real^{q}$ and a matrix $\H \in \Real^{q \times n}$.
A more general SOCP formulation often includes a Cartesian product of second-order cones.
However, only one second-order cone is enough for the discussions in this paper.

Yamashita et al.~\cite{yamashita2015efficient} introduced a new vector $\z := \A \x \in \Real^Z$,
and converted the quadratic constraint $\x^T \A \x \le 2 \theta$ 
into $||\B \z|| \le \sqrt{2 \theta}$ with a matrix $\B \in \Real^{Z \times Z}$ that satisfies $\B^T \B = \A^{-1}$.
Though the Wright numerator  matrix $\A$ itself is not a sparse matrix,
the matrices $\A^{-1}$ and $\B$ possess favorable sparsity.
The computation time reduction reported in~\cite{yamashita2015efficient} was mainly derived from these sparsity.
Using these new vector $\z$ and matrix $\B$, we transformed the ED problem (\ref{eq:ED})
into the following SOCP problem with integer constraints;
\begin{eqnarray*}
%(SOCP) : 
\begin{array}{lcl}
\max &:& (\A^{-1}\g)^T \z \\
\mbox{subject to} &:& (\A^{-1}\e)^T \z = 1,\\
& & \left(\begin{array}{c} \sqrt{2 \theta} \\ \B \z \end{array}\right) \in \KC^{1+Z},\\
& & \begin{array}{lcl}
[\A^{-1} \z]_i \in \left\{0, \frac{1}{N}\right\} & \mbox{for} & i \in V,\\ 
 {}   [\A^{-1} \z]_i = c_i & \mbox{for} & i \in F. 
\end{array}
\end{array}\label{eq:MI-ED}
\end{eqnarray*}
Here, we use the notation $[\A^{-1} \z]_i$ to denote the $i$th element of $\A^{-1} \z$.
It may seem that we would remove $\x_F$ from this formulation by fixing $\x_F = \c_F$
and reduce the sizes the problem. % we need to solve.
However, such elimination would strongly diminish the efficiency of the SOCP problem,
since it completely destroys the favorable 
sparsity that appear in $\A^{-1}$ and $\B$.

By applying the continuous relaxation to the binary constraints, 
we obtain an SOCP relaxation problem of the ED problem;
\begin{eqnarray}
\begin{array}{lclcl}
OPT_{SOCP} &:=& \max &:& (\A^{-1}\g)^T \z\\
& & \mbox{subject to} &:& (\A^{-1}\e)^T \z = 1,\\
& & & & \left(\begin{array}{c} \sqrt{2 \theta} \\ \B \z \end{array}\right) \in \KC^{1+Z},\\
& & & & \begin{array}{lcl}
0 \le [\A^{-1} \z]_i \le \frac{1}{N} & \mbox{for} & i \in V,\\ 
 {}   [\A^{-1} \z]_i = c_i & \mbox{for} & i \in F. 
\end{array}
\end{array}\label{eq:SOCP-R}
\end{eqnarray}

% The optimal values of  ED~(\ref{eq:ED}), $SDP_R$~(\ref{eq:SDP-R}), $LP_R$~(\ref{eq:LP-R})
% and $SOCP_R$~(\ref{eq:SOCP-R}) will be denoted by
$OPT_{ED}$, $OPT_{SDP}$, $OPT_{LP}$ and $OPT_{SOCP}$, respectively.
From the derivation of the LP relaxation problem~(\ref{eq:LP-R}),
it is natural that 
the SDP relaxation problem~(\ref{eq:SDP-R}) gives closer an optimal value 
than the LP relaxation problem, that is,  we know $OPT_{SDP} \le OPT_{LP}$.
In contrast, the relation of the SOCP relaxation~(\ref{eq:SOCP-R}) is not so explicit,
since the SOCP relaxation was derived 
by a continuous relaxation independently from the SDP or LP relaxation.

The strength of these relaxation problems can be summarized in Lemma~\ref{le:relaxation}.
For the discussion there, we prepare some notation and introduce an assumption.
We use $\SC_m(\v)$ to denote the sum of the $m$ smallest elements of $\v \in \Real^n$.
More precisely, when $\hat{v}_{1} \le \hat{v}_2 \le \ldots \le \hat{v}_n$ is the sorted vector of $\v$ in the ascending order,
the definition of $\SC_m(\v)$ is given by  $\SC_m(\v) := \sum_{i=1}^m \hat{v}_i$. 
The symbol $\hat{A}_{\hat{W}}$ indicates the set of the collection of $\A_{VV}$ with respect to $\hat{W}$, 
that is, $\hat{A}_{\hat{W}} := \left\{A_{ij} : (i,j) \in \hat{W}\right\}$.
We define a vector $\hat{\y}_V \in \Real^{|V|}$ by $[\hat{y}_V]_i := 1$ for $i = 1, \ldots, N-p$
and $[\hat{y}_V]_i := -1$ for $i = N-p+1, \ldots, |V|$. 
This vector satisfies $\e_V^T \hat{\y}_V = \bar{N}$.
In the following this discussion, we make the following assumption on the input data 
of the ED problem~(\ref{eq:ED}).
From preliminary numerical tests, 
we verified that this assumption holds for practical datasets of pine orchards 
and datasets generated by simulations. The details of these dataset will be described
in Section~5.
\begin{ASSUM}\label{as:A}
The input data of (\ref{eq:ED}) satisfies 
\begin{eqnarray*}
\SC_{\hat{N}}(\hat{A}_{\hat{W}}) \le \frac{2 \bar{\theta} - 2\mbox{Trace}(\A_{VV}) + \e_V^T \A_{VV} \e_V - 2 \bar{\c}_F^T \hat{\y}}{4},
\end{eqnarray*}
where $\hat{N} := \frac{\bar{N}^2 + |V|^2 - 2|V|}{4}$.
\end{ASSUM}
We should ensure that $\hat{N}$ is a positive integer, otherwise we need to manage 
a fractional number in the definition of $\SC$.
The positiveness is derived from
$\bar{N}^2 + |V|^2 - 2|V| \ge \bar{N}^2 + 1 \ge 1$
by $|V| \ge 2$ of Remark~\ref{re:remarkPN},
and $\hat{N}$ is integer by
\begin{eqnarray*}
\bar{N}^2 + |V|^2 - 2|V| &=& \left\{2N (1- \e_F^T \c_F) - |V| \right\}^2 + |V|^2 - 2|V| \\
&=& \left\{2N (1- \frac{p}{N}) - |V| \right\}^2 + |V|^2 - 2|V| \\
&=& 4 \left\{(N-p)^2 -  |V| (N-p) +  \frac{|V| (|V| -1)}{2}\right\}.
\end{eqnarray*}

We are now prepared to examine the relation between the relaxation problems.
\begin{LEMM}\label{le:relaxation}
It holds for the optimal values of the relaxation problems that 
\begin{eqnarray*}
OPT_{ED} \le OPT_{SDP} \le OPT_{SOCP}. % \le OPT_{LP}.
\end{eqnarray*}
Furthermore, if Assumption~\ref{as:A} holds, then
\begin{eqnarray*}
OPT_{ED} \le OPT_{SDP} \le OPT_{SOCP}  \le OPT_{LP}.
\end{eqnarray*}
\end{LEMM}

\noindent \textbf{Proof:} 
[$OPT_{ED} \le OPT_{SDP}$]
When we derived (\ref{eq:SDP-R}), we ignored the rank-1 constraint in (\ref{eq:ED4}).
From this derivation, for any feasible solution $\x \in \Real^Z$ of (\ref{eq:ED}),
the corresponding vector $\y_V \in \Real^{|V|}$ 
through the connections $\x = \left(\begin{array}{cc} \x_V \\ \c_F \end{array}\right)$, then
$\y_V = 2 N \x_{V} - \e_{V}$
is also a feasible solution of (\ref{eq:SDP-R}).
Furthermore, from these connections hold, it holds that 
$\g^T \x = \bar{\g}^T \y_V + \bar{g}$.
The objective functions 
of (\ref{eq:ED4}) and (\ref{eq:SDP-R}) are same and
the feasible region of (\ref{eq:SDP-R}) is wider than that of 
(\ref{eq:ED4}) substantially, hence, we have $OPT_{ED} \le OPT_{SDP}$. 

[$OPT_{SDP} \le OPT_{SOCP}$] We take any feasible solution 
$\y_V \in \Real^{|V|}$ and $\Y_{VV} \in \SMAT^{|V|}$ of (\ref{eq:SDP-R}).
It is enough to check that 
$\z = \A \left(\begin{array}{c} \frac{\y_V + \e_V}{2N} \\ \c_F \end{array}\right)$
is a feasible solution of (\ref{eq:SOCP-R}).

From 
$\left(\begin{array}{cc} -2\bar{N} & \e_V^T \\ \e_V & \O \end{array}\right)
\bullet
\left(\begin{array}{cc} 1 & \y_V^T \\ \y_V & \Y_{VV} \end{array}\right)
= 0$, we obtain $\e_V^T \y_V = \bar{N} = 2N (1-\e_F^T \c_F) - |V|$, hence, 
\begin{eqnarray*}
(\A^{-1} \e)^T \z  =   
\left(\begin{array}{c}\e_V \\ \e_F \end{array}\right)^T
\left(\begin{array}{c} \frac{\y_V + \e_V}{2N} \\ \c_F \end{array}\right)
= \frac{\e_V^T \y_V + |V|}{2N} + \e_F^T \c_F = 1.
\end{eqnarray*}
By applying the Schur complement to  the positive semidefinite condition
$\left(\begin{array}{cc} 1 & \y_V^T \\ \y_V & \Y_{VV} \end{array}\right)
\succeq \O$, it holds $\Y_{VV} - \y_V \y_V^T \succeq \O$.
Since $\A \bullet \X \ge 0$ holds for any two positive semidefinite matrices of the same dimension 
$\A$ and $\X$~\cite{TODD01} and the Wright numerator matrix is always positive definite,
it holds $\A_{VV} \bullet (\Y_{VV} - \y_V \y_V^T) \ge 0$, therefore,
$\A_{VV} \bullet \Y_{VV}  \ge \y_V^T \A_{VV} \y_V$.
Using the relation
$\left(\begin{array}{cc} -2\bar{\theta} & \bar{\c}_F^T \\ \bar{\c}_F & \A_{VV} \end{array}\right)
\bullet
\left(\begin{array}{cc} 1 & \y_V^T \\ \y_V & \Y_{VV} \end{array}\right)
\le 0$, 
we obtain $\y_V^T \A_{VV} \y_V + 2 \bar{\c}_F^T \y_V \le 2 \bar{\theta}$.
From the definitions of $\y_V, \bar{\c}_F, \bar{\theta}, \B$ and $\z$,
we can derive $\z \B^T \B \z \le 2 \theta$, therefore, 
$\left(\begin{array}{c} \sqrt{2 \theta} \\ \B \z \end{array}\right) \in \KC^{1+Z}$. 
From $\Y_{VV} - \y_V \y_V^T \succeq \O$, we also have
$Y_{ii} \ge y_i^2$  for $ i = 1, \ldots, |V|$.
Furthermore, due to the constraint 
$\left(\begin{array}{cc} -1 & \0^T \\ \0 & \e_i \e_i^T  \end{array}\right)
\bullet
\left(\begin{array}{cc} 1 & \y_V^T \\ \y_V & \Y_{VV} \end{array}\right)
= 0$, it holds $Y_{ii} = 1$ for $i = 1, \ldots, |V|$,
 consequently $-\e_V \le \y_V \le \e_V$.
From $\A^{-1} \z = \left(\begin{array}{c} \frac{\y_V + \e_V}{2N} \\ \c_F \end{array}\right)$, it is now clear that 
$0 \le [\A^{-1} \z]_i \le \frac{1}{N}$  for  $i \in V$
and that  $[\A^{-1} \z]_i = c_i$ for $i \in F$.
Furthermore, the objective value of (\ref{eq:SDP-R}) at $\y_V$ is same as 
that of (\ref{eq:SOCP-R}) at $\z$ if
$\z = \A \left(\begin{array}{c} \frac{\y_V + \e_V}{2N} \\ \c_F \end{array}\right)$.
Hence, we obtain $OPT_{SDP} \le OPT_{SOCP}$.

[$OPT_{SOCP} \le OPT_{LP}$] 
%We now take any feasible solution $\z \in \Real^Z$ of $SOCP_R$~(\ref{eq:SOCP-R}).
%This time, we first introduce $\x = \A^{-1} \z$.
%Since $\z$ is a feasible solution of $SOCP_R$, it holds that
%$\e^T \x = 1$, $\x^T \A \x \le 2 \theta$,
%$0 \le \x_V \le \frac{1}{N}$ and $\x_F = \c_F$.
%When we define a vector $\y_V := 2 N \x_V - \e_V$ and 
%a matrix $\Y_V \in \SMAT^{|V|}$ such that $Y_{ii} = 1$ for $i = 1, \ldots, |V|$
%and $Y_{ij} = Y_{ji} = y_i y_j$ for $i = 1, \ldots |V|$ and $j = i+1, \ldots, |V|$.
We first consider an LP problem
\begin{eqnarray}
\begin{array}{lcl}
\min &:& \c^T \futoji{\eta}\\
\mbox{subject to} &:& \sum_{i=1}^n \eta_i = K,\\
& & 
0 \le \eta_i \le 1 \mbox{ for } i = 1, \ldots, n,
\end{array}\label{eq:simple-lp}
\end{eqnarray}
where the decision variable is $\futoji{\eta} \in \Real^n$ and the input vector is $\c\in \Real^n$ and $K$ is a positive integer.
The optimal value of this LP problem is $\SC_K(\c)$ and this value can be attained 
at $\hat{\futoji{\eta}} \in \Real^n$ such that $\hat{\eta}_i = 1$ for $i = 1, \ldots, K$ and $\hat{\eta}_i = 0$ for $i = K+1, \ldots, n$.

If we ignore the quadratic constraint of (\ref{eq:SOCP-R}) and we reverse the variable into $\x = \A^{-1} \z$, 
we obtain an optimization problem of form
\begin{eqnarray}
% (SOCP_R) : 
\begin{array}{lcl}
\max &:& \g_V^T \x_V + \g_F^T \c_F,\\
\mbox{subject to} &:& \e_V^T \x_V  = 1 - \frac{p}{N},\\
& & \begin{array}{lcl}
0 \le x_i \le \frac{1}{N} & \mbox{for} & i \in V.
\end{array}
\end{array}\label{eq:SOCP-R2}
\end{eqnarray}
Since $g_1 \ge g_2 \ge \ldots g_{|V|}$, 
the optimal value of (\ref{eq:SOCP-R2}) is given by $-\frac{\SC_{N-p}(-\g_V)}{N} + \g_F^T \c_F$
 in a similar way to (\ref{eq:simple-lp})
and  an optimal solution is $\hat{\x}_V := \frac{\hat{\y}_V + \e_V}{2N}$.
Therefore, it holds that $OPT_{SOCP} \le \g_V^T \hat{\x}_V + \g_F^T \c_F$.

% To show the existence of $\hat{\Y}_{VV} \in \SMAT^{|V|}$ such that 
% the pair of $\hat{\y}_V$ and $\hat{\Y}_{VV}$ is a feasible solution of $(LP_R)$, 
Next, we define $\rho_{LP}$ to denote
the optimal value of the following LP problem;
\begin{eqnarray}
\begin{array}{lclcl}
\rho_{LP} &:=&  \min &:& \A_{VV} \bullet \Y_{VV} \\
& & \mbox{subject to} &:& (\e_V \e_V^T) \bullet \Y_{VV}  = \bar{N}^2, \\
& & & & Y_{ii} = 1 \mbox{ for } i \in V,\\
& & & & -1 \le Y_{ij} \le 1 \mbox{ for } (i,j) \in \hat{W},\\
& &  & & \Y_{VV} \in \SMAT^{|V|}.
\end{array}\label{eq:LP-Feas}
\end{eqnarray}
We convert this problem introducing $\bar{X}_{ij} := \frac{Y_{ij} + 1}{2}$ 
for $(i,j) \in \hat{W}$.
The following LP problem is equivalent to (\ref{eq:LP-Feas}), therefore,
its optimal value must be $\rho_{LP}$.
\begin{eqnarray}
\begin{array}{lclcl}
\rho_{LP} &=& \min &:& 4 \sum_{(i,j) \in \hat{W}} A_{ij} \bar{X}_{ij} 
- 2 \sum_{(i,j) \in \hat{W}} A_{ij} + \mbox{Trace}(\A_{VV})\\
& & \mbox{subject to} &:& \sum_{(i,j) \in \hat{W}}\bar{X}_{ij} = \hat{N},\\
& & & & 0 \le \bar{X}_{ij} \le 1 \mbox{ for } (i,j) \in \hat{W}. \\
\end{array}\label{eq:LP-Feas2}
\end{eqnarray}
The structure of this problem is same as (\ref{eq:simple-lp}), hence, 
it holds that
$\rho_{LP} = 4 \SC_{\hat{N}}(\hat{A}_{\hat{N}}) 
- \e_V^T \A_{VV}  \e_V+ 2 \mbox{Trace}(\A_{VV})$.

Let $\hat{\Y}_{VV}$ be a part of an optimal solution of (\ref{eq:LP-Feas}).
From Assumption~\ref{as:A}, it holds that
\begin{eqnarray*} 
\A_{VV} \bullet \hat{\Y}_{VV} + 2 \bar{\c}_F^T \hat{\y}_V 
= \rho_{LP} + 2 \bar{\c}_F^T \hat{\y}_V  
= 4 \SC_{\hat{N}}(\hat{A}_{\hat{N}}) 
- \e_V^T \A_{VV}  \e_V+ 2 \mbox{Trace}(\A) + 2 \bar{\c}_F^T \hat{\y}_V 
\le 2 \bar{\theta}.
\end{eqnarray*}
Furthermore, $\hat{\y}_V$ satisfies $-1 \le \hat{y}_i \le 1$ for $i \in V$ and $\e_V^T \hat{\y}_V = \bar{N}$ by its definition and 
$\hat{\Y}_{VV}$ satisfies all the constraints of (\ref{eq:LP-Feas}).
Consequently, the pair $\hat{\y}_V$ and $\hat{\Y}_{VV}$
is a feasible solution of (\ref{eq:LP-R}) 
and this leads to the inequality we wanted to obtain.
\begin{eqnarray*}
OPT_{LP} \ge \left(\begin{array}{cc} 0 & \bar{\g}_V^T \\ \bar{\g}_V & \O \end{array}\right)
\bullet
\left(\begin{array}{cc} 1 & \hat{\y}_V^T \\ \hat{\y}_V & \hat{\Y}_{VV} \end{array}\right)
+ \bar{g} = 2 \bar{\g}_V ^T \hat{\y}_V + \bar{g} = \g_V^T \hat{\x}_V + \g_F^T \c_F \ge OPT_{SOCP}.
\end{eqnarray*}
\qed

\begin{REMA}\label{re:LP}
Since an optimal solution of a further relaxation problem of (\ref{eq:LP-R}) 
\begin{eqnarray*}
% (SOCP_R) : 
\begin{array}{lcl}
\max &:& 2 \bar{\g}_V^T \y_V + \bar{g} \\
\mbox{subject to} &:& \e_V^T \y_V  = \bar{N}, \\
& & -1 \le y_i \le 1  \mbox{ for }  i \in V
\end{array}\label{eq:LP-R2}
\end{eqnarray*}
is $\hat{\y}_V$, its optimal value 
$2 \bar{\g}_V^T \hat{\y}_V + \bar{\g}$ must be an upper bound
of $OPT_{LP}$. 
On the other hand,  from the proof of Lemma~\ref{le:relaxation},
when Assumption~\ref{as:A} holds, there exists some
$\hat{\Y}_{VV} \in \SMAT^{|V|}$ 
such that the pair of $\hat{\y}_V$ and $\hat{\Y}_{VV}$ is a feasible 
solution of (\ref{eq:LP-R}) with the objective value 
$2 \bar{\g}_V^T \hat{\y}_V + \bar{\g}$. 
Therefore, $\hat{\y}_V$ is also an optimal solution of (\ref{eq:LP-R}).
This indicates that we can obtain the solution of 
(\ref{eq:LP-R}) at the computation cost for sorting $\g_V$ 
instead of solving (\ref{eq:LP-R}) as an LP problem, when Assumption~\ref{as:A} holds.
\end{REMA}
This remark implies that the LP relaxation (\ref{eq:LP-R}) is not so tight 
against the original ED problem (\ref{eq:ED}).
In contrast, we observed through preliminary numerical tests that 
the vector $\left(\begin{array}{c} \hat{\x}_V \\ \c_F \end{array}\right)$ defined with
an optimal solution $\hat{\x}_V$  of (\ref{eq:SOCP-R2}),
is not always a feasible solution of (\ref{eq:SOCP-R}) even if Assumption~\ref{as:A} holds.
Therefore, the feasible region of the SOCP relaxation problem is strictly narrower than that of the LP relaxation problem,
and the SOCP relaxation gives a tighter approximation than the LP relaxation in general,
even though the relaxation were derived independently.

\begin{REMA}\label{re:SDP}
The SDP relaxation problem (\ref{eq:SDP-R}) has no interior-feasible point.
\end{REMA}
If a pair  $\y_V \in \Real^{|V|}$ and $\Y_{VV} \in \SMAT^{|V|}$ satisfies all the constraint of (\ref{eq:SDP-R}),
the pair is a feasible point. When 
the matrix $\left(\begin{array}{cc} 1 & \y_V^T \\ \y_V & \Y_{VV} \end{array}\right)$
is a positive definite matrix for some feasible point  $\y_V \in \Real^{|V|}$ and $\Y_{VV} \in \SMAT^{|V|}$,
we say that (\ref{eq:SDP-R}) has an interior-feasible point.
We can show that $\left(\begin{array}{cc} 1 & \y_V^T \\ \y_V & \Y_{VV} \end{array}\right)$ is not positive definite 
for any feasible point of (\ref{eq:SDP-R}).
To show this, we take a feasible point 
$\y_V \in \Real^{|V|}$ and $\Y_{VV} \in \SMAT^{|V|}$. Then, we have $\e_V^T \y_V = \bar{N}$ and
$(\e_V \e_V^T) \bullet \Y_{VV} = \bar{N}^2$.
If $\bar{N} \ne 0$, it holds 
\begin{eqnarray*}
\left(\begin{array}{c} 1  \\ -\e_V/\bar{N}  \end{array}\right)^T 
\left(\begin{array}{cc} 1 & \y_V^T \\ \y_V & \Y_{VV} \end{array}\right)
\left(\begin{array}{c} 1 \\  -\e_V/\bar{N}  \end{array}\right)
= 1 - 2 \e_V^T \y_V / \bar{N} + \e_V^T \Y_{VV} \e_V / \bar{N}^2 = 0.
\end{eqnarray*}
In addition, for the case $\bar{N} = 0$, it holds 
\begin{eqnarray*}
\left(\begin{array}{c} 0  \\ \e_V  \end{array}\right)^T 
\left(\begin{array}{cc} 1 & \y_V^T \\ \y_V & \Y_{VV} \end{array}\right)
\left(\begin{array}{c} 0 \\  \bar{N}  \end{array}\right)
= \e_V^T \Y_{VV} \e_V = \bar{N}^2 = 0.
\end{eqnarray*}
In either case, there exists a nonzero vector that makes the quadratic from zero,
the matrix is not positive definite, therefore, 
(\ref{eq:SDP-R}) has no interior-feasible point.

\section{Theoretical evaluation of the SDP relaxation problems with a randomized algorithm}
The solutions obtained from the conic relaxation problem~(\ref{eq:SDP-R}),  (\ref{eq:LP-R}) and (\ref{eq:SOCP-R})
are not always a feasible solution of the ED problem~(\ref{eq:ED}), 
since we ignored some constraints of the NP-hard problem 
to derive the conic relaxation problems that are solvable in polynomial time.
When SDP relaxation approaches are used,  examine randomized algorithms often follow 
to generate feasible solutions.
A randomized algorithm using the solutions obtained through SDP relaxation problems
was first introduced for max-cut problems
in~\cite{GOEMANS95}. They showed that the expectation objective value 
obtained by their randomized algorithm on average was at least 0.878 of that of 
an SDP relaxation problem. 
Since the optimal value of a max cut problem exists between an objective value
of any feasible solution and the value obtained from the SDP relaxation problem,
their algorithm have an expected approximation
factor of 0.878.
Many researches followed~\cite{GOEMANS95} to extend its results to more general 
quadratic-constraint problems using the framework of SDP relaxation methods.
Among them, Tseng~\cite{tseng2003further} discussed one of the most general cases
and gave its probabilistic analysis.
Wu et al.~\cite{wu2013new} also analyzed the expectation values using a different
randomized algorithm.

In this section, 
we employ the result of~\cite{tseng2003further} to give
theoretical bounds on 
the expected objective value of a randomized algorithm. 
% based on the SDP relaxation problems~(\ref{eq:SDP-R}) 
% is $1 - \frac{\pi}{2}$ up to constant numbers.
Tseng~\cite{tseng2003further} applied the SDP relaxation methods
to a quadratically-constrained quadratic programming (QCQP) problem:
\begin{eqnarray}
\begin{array}{lclcl}
\overline{OPT}_{QCQP} &:= &  \max &:& \y^T \A^0 \y + (\b^0)^T \y + c^0 \\
&  & \mbox{subject to} &:& \y^T \A^k \y + (\b^k)^T \y + c^k \le 0 \mbox{ for } k=1,\ldots, m.
\end{array}
\label{eq:QCQP}
\end{eqnarray}
Here, the variable is $\y \in \Real^n$, while the input data are 
$\A^0, \ldots, \A^m \in \SMAT^n$, $\b^0, \ldots, \b^m \in \Real^n$
and $c^0, c^1 \ldots, c^m \in \Real$. For simplicity, the constant in the objective function is fixed to $c^0 = 0$.
The  QCQP originally discussed in \cite{tseng2003further} is a minimization problem,  but
we consider a maximization problem since the ED problem (\ref{eq:ED}) is a maximization problem.

% We now use $\overline{OPT}_{QCQP}$ to denote the optimal value of 
When we apply the lift-and-project method of Lov{\'a}sz and 
Schrijver~\cite{lovasz1991cones}
to (\ref{eq:QCQP}), the resultant SDP relaxation problem
is given as follows:
\begin{eqnarray}
\begin{array}{lclcl}
\overline{OPT}_{SDP} &:=&  \max &:& \B^0 \bullet \Y \\
& &  \mbox{subject to} &:& \B^k \bullet \Y \le 0  \mbox{ for } k=1,\ldots, m,\\
& & & & \B^{m+1} \bullet \Y = 1 , \quad \Y \succeq \O
\end{array}
\label{eq:QCQP-SDP}
\end{eqnarray}
where $\B^k := \left(\begin{array}{cc}
c^k & (\b^k)^T \\ \b^k & \A^k 
\end{array}\right)$ for $k = 0, \ldots, m$ and
$\B^{m+1} := \left(\begin{array}{cc}
1 & \0^T \\ \0 & \O
\end{array}\right)$,
and 
the decision variable is $\Y \in \SMAT^{1+n}$.
In the following discussions, we start the row or column index of $\Y$ from zero,
therefore, the elements of $\Y$ are denoted by $Y_{00}, Y_{01}, \ldots, Y_{nn}$.
% We denote, by  $\overline{OPT}_{SDP}$, the optimal value of (\ref{eq:QCQP-SDP}). 
It is known that if we add the rank-1 constraint $\mbox{rank}(\Y) = 1$ to 
(\ref{eq:QCQP-SDP}), 
the two problems (\ref{eq:QCQP}) and (\ref{eq:QCQP-SDP})
are equivalent.
In other words, we ignored the rank-1 constraint from (\ref{eq:QCQP}) 
to derive (\ref{eq:QCQP-SDP}), hence,
$\overline{OPT}_{QCQP} \le \overline{OPT}_{SDP}$.

The randomized algorithm of~\cite{tseng2003further} can be summarized 
as follow.
We assume that (\ref{eq:QCQP}) and (\ref{eq:QCQP-SDP})
are feasible and that (\ref{eq:QCQP-SDP}) has an optimal solution, denoted as $\Y^*$.
This solution $\Y^*$ is factorzied 
with a matrix 
$\V \in \Real^{(1+n) \times (1+n)}$ such that  $\Y^* = \V^T \V$. 
Such $\V$ is available,  for example, by
the Cholesky factorization or the eigenvalue decomposition.
We use $\v^0, \v^1, \ldots, \v^{n} \in \Real^{1+n}$ to denote the columns of $\V$.
Then, a vector $\v \in \Real^{1+n}$ is chosen randomly 
from the unit sphere in $\Real^{1+n}$ based on uniform distribution.
Finally, the randomized algorithm outputs a solution $\tilde{\y} \in \Real^{1+n}$ 
defined by
\begin{eqnarray*}
 \tilde{y}_i := \sqrt{Y^*_{ii}} \mbox{sign}(\v^T \v^0) \mbox{sign}(\v^T \v^i)
\mbox{ for } i = 0, \ldots, n
\end{eqnarray*}
where $\mbox{sign}(a) = 1$ if $a \ge 0$ and $\mbox{sign}(a) = -1$ if $a < 0$.
We remark that from the definition of $\B^{m+1}$, 
it always holds that $Y_{00}^* = 1$, hence,
$\tilde{y}_0 = \sqrt{1} (\mbox{sign}(\v^T \v^0))^2 = 1$.

The set $\IC$ is introduced to indicate diagonal-matrix constraints of (\ref{eq:QCQP});
\begin{eqnarray*}
\IC := \left\{ k \in \{1, 2, \ldots, m\} :
\A^k \mbox{ is a diagonal matrix and } \b^k = \0 \right\}.
\end{eqnarray*}
To measure 
% the violation on the $l$th constraint of (\ref{eq:QCQP-SDP}) for $l \notin \IC$ or 
a  shift in the objective function, % with $l=0$ , 
$\rho_{SDP}^0$ is defined as the optimal value of 
the following SDP problem
\begin{eqnarray}
\begin{array}{lclcl}
 \rho_{SDP}^0 &:=& \min &:& \B^0 \bullet \Y \\
& &  \mbox{subject to} &:& \B^k \bullet \Y =\B^k \bullet \Y^*
\mbox{ for } k \in \IC,  \\
& & & & \B^{m+1} \bullet \Y = 1 , \quad \Y \succeq \O.
\end{array}\label{eq:SDP-bound}
\end{eqnarray}

Tseng~\cite{tseng2003further} showed a relation
between the expected objective value 
of the generated solution $\tilde{\y}$ and 
the optimal values of the SDP  problems.
\begin{THEO}\label{th:Tseng}
% $\mathrm{[\cite{tseng2003further}, Theorem2]}$
{\upshape\cite[Theorem 2]{tseng2003further}}
If the SDP relaxation problem~(\ref{eq:QCQP-SDP})
has an optimal solution $\Y^*$ and  a set 
$\left\{\y \in \Real^n : \y^T \A^k \y + (\b^k)^T \y + c^k \le 0, k \in \IC\right\}$ 
is bounded, then
\begin{eqnarray*}
E[ \tilde{\y}^T \A^0 \tilde{\y} + (\b^0)^T \y]
\ge \frac{2}{\pi} \overline{OPT}_{SDP} 
+ \left(1-\frac{2}{\pi}\right) \rho_{SDP}^0.
\end{eqnarray*}
\end{THEO}

Let us return to the ED problem~(\ref{eq:ED}).
We analyze the performance 
of the output solution $\tilde{\y}_V$ that is generated 
by the above randomized algorithm using 
the optimal solution $\Y^*$  of 
the SDP relaxation problem~(\ref{eq:SDP-R}).
From the form of (\ref{eq:SDP-R}), 
the objective value at $\tilde{\y}_V$  is $2 \g_V^T \tilde{\y} + \bar{g}$.
The following lemma provides a theoretical aspects on the expected value
of this objective function.

\begin{LEMM}\label{le:SDP-R}
For the ED problem (\ref{eq:ED}), 
the expected objective value obtained through the randomized algorithm  
is bounded by 
\begin{eqnarray*}
\frac{2}{\pi} OPT_{SDP}
+ \left( 1- \frac{2}{\pi} \right) (-2 \bar{\g}_V^T \e_V  + \bar{g})
 \le E[2\bar{\g}_V^T \tilde{\y}_V + \bar{g}]  \le 
\alpha OPT_{SDP} +  (1-\alpha) (2 \bar{\g}_V^T \e_V + \bar{g}),
\end{eqnarray*}
where $\alpha := \min\left\{\frac{2}{\pi} \frac{\theta}{1-\cos \theta}
: 0 \le \theta \le \pi \right\} \approx 0.878$.
\end{LEMM}

\textbf{Proof:}

First, we derive the lower bound of the objective function by use of 
 Theorem~\ref{th:Tseng}.
To embed the SDP relaxation problem~(\ref{eq:SDP-R}) arising from the ED problem 
into the framework developed in~\cite{tseng2003further}, 
we embed the variable vector $\y_V$ and matrix $\Y_{VV}$ into the matrix $\Y \in \SMAT^{1+|V|}$ as
$%\begin{eqnarray*}
\Y = \left(\begin{array}{cc} 
Y_{00} & \y_V^T \\ \y_V & \Y_{VV}
\end{array}\right).$ %\end{eqnarray*}
In particular, we identify $y_i = Y_{0i} = Y_{i0}$ for $i = 1, \ldots, |V|$.
%For the consistent notation, we index the row or column of $\Y$ as $\{0, 1, 2, \ldots, |V|\}$.
For the input matrices $\B^0, \ldots, \B^{2|V|+5}$,
we prepare
\begin{eqnarray*}
\left\{
\begin{array}{lllll}
c^0 = 0, & \b^0 = \bar{\g}_V, & \A^0 = \O, \\
c^1 = -2 \bar{\theta}, & \b^1 = \bar{\c}_F, & \A^1 = \A_{VV}, \\
c^2 = -2\bar{N}, & \b^2 = \e_V, & \A^2 = \O, \\
c^3 = 2\bar{N}, & \b^3 = -\e_V, & \A^3 = \O, \\
c^4 = -\bar{N^2}, & \b^4 = \0, & \A^4 = \e_V \e_V^T, \\
c^5 = \bar{N^2}, & \b^5 = \0, & \A^5 = -\e_V \e_V^T, \\
c^{5+i} = -1, & \b^{5+i} = \0, & \A^{5+i} = \e_i \e_i^T  \mbox{ for }  i \in V,\\
c^{5+|V| + i} = 1, & \b^{5+|V| + i} = \0, & \A^{5+|V|+i} = -\e_i \e_i^T  \mbox{ for }  i \in V. \\
\end{array}
\right.
\end{eqnarray*}

The number of input matrices in the form of (\ref{eq:QCQP-SDP}) is $m = 2|V|+5$.
For example, $\B^{5+i} \bullet \Y \le 0$ and $\B^{5+|V|+i} \bullet \Y \le 0$ lead to 
$Y_{ii} = 1$ for $i \in V$. In addition,  $Y_{00} = 1$ is guaranteed by $\B^{m+1} \bullet \Y = 1$.

The set of diagonal constraints is $\IC =  \{ 5+i : i \in V\} \cup
\{ 5 + |V| + i : i \in V\}$.
From this $\IC$, the feasible set of (\ref{eq:SDP-bound}) is given by
$\FC := \left\{\Y \in \SMAT^{1+|V|} : Y_{ii} = 1 \mbox{ for } i = 0, 1, \ldots, |V|
\mbox{ and } \Y \succeq \O\right\}$.
Hence, we obtain
\begin{eqnarray*}
\rho_{SDP}^0 = \min \left\{ \B^0 \bullet \Y : \Y \in \FC \right\} 
 = \min \left\{ 2 \sum_{i \in V}  \bar{g}_i Y_{0,i} : \Y \in \FC \right\} 
% = - 2 \sum_{i \in V} \bar{g}_i 
= - 2 \bar{g}_V^T \e_V.
\end{eqnarray*}
Here, 
a combination of 
the matrix-completion method~\cite{FUKUDA00, NAKATA03, yamashita2015fast}
with a property $\bar{\g}_V \ge \0$ ensures that 
an optimal solution of this minimization problem is given as
$\Y = \left(\begin{array}{cc} 1 & -\e_V^T \\ 
-\e_V & \e_V \e_V^T \end{array}\right)$.

We should note that the SDP relaxation problem~(\ref{eq:SDP-R}) has a constant term
$\bar{g}$ in the objective function, but we have to set $c^0 = 0$ to employ
Theorem~\ref{th:Tseng}.
By taking the shift of $\bar{g}$ into account, Theorem~\ref{th:Tseng} gives  a lower bound;
\begin{eqnarray*}
E[2\bar{\g}_V^T \tilde{\y}_V + \bar{g}] &=&  
E[2\bar{\g}_V^T \tilde{\y}_V] + \bar{g} \\
&\ge &
\frac{2}{\pi} (\overline{OPT}_{SDP}) 
 -  \left( 1- \frac{2}{\pi} \right) 2 \bar{\g}_V^T \e_V + \bar{g}\\
&=&
\frac{2}{\pi} (OPT_{SDP} - \bar{g}) 
-  \left( 1- \frac{2}{\pi} \right) 2 \bar{\g}_V^T \e_V + \bar{g}\\
&=& \frac{2}{\pi} OPT_{SDP}
+ \left( 1- \frac{2}{\pi} \right) (-2 \bar{\g}_V^T \e_V  + \bar{g}).
\end{eqnarray*}

To consider an upper bound, we first evaluate $E[\tilde{y}_i]$  for $i \in  V$.
From $Y_{ii}^*= 1$ for $i \in \{0\} \cup V$ in (\ref{eq:SDP-R}) 
and  the definition of $\tilde{\y}_V$,
and it holds that 
$\tilde{y}_i = 1$ if $\mbox{sign}(\v^T \v^0) = \mbox{sign}(\v^T \v^i)$,
and $\tilde{y}_i = -1$ if $\mbox{sign}(\v^T \v^0) = -\mbox{sign}(\v^T \v^i)$.
The discussion in \cite{GOEMANS95} indicates that the probability 
of the event $\mbox{sign}(\v^T \v^0) = \mbox{sign}(\v^T \v^0)$
is given as $1 - \frac{1}{\pi} \mbox{arccos} (Y_{0i}^*)$.
Therefore,  we have
\begin{eqnarray*}
E[\tilde{y}_i] &= &
1 \cdot \left(1 - \frac{1}{\pi} \mbox{arccos} (Y_{0i}^*)\right)
+ (-1) \cdot \left\{1-\left(1 - \frac{1}{\pi} \mbox{arccos} (Y_{0i}^*)\right)\right\} \\
&=& 1-\frac{2}{\pi} \mbox{arccos} (Y_{0i}^*) 
\le  \alpha (Y_{0i}^* - 1) + 1.
\end{eqnarray*}
The last inequality was derived from 
the
inequality 
$\frac{\mbox{arccos}(y)}{\pi} \ge \alpha \frac{1-y}{2}$
 for $-1 \le y \le 1$
(Lemma~3.4 of~\cite{GOEMANS95})
and
$-1 \le Y_{0i}^* \le 1$ 
due to $\Y^* \succeq \O$ and $Y_{00}^* = Y_{ii}^* = 1$.

As a result, we obtain an inequality
\begin{eqnarray*}
E[2\bar{\g}_V^T \tilde{\y}_V]  + \bar{g}
&=& 2 \sum_{i \in V} \bar{g}_i E[\tilde{y}_i]  + \bar{g} \\
&\le& 2 \alpha \sum_{i \in V} \bar{g}_i  Y_{0i}^* + 2 (1-\alpha) \g_V^T \e_V 
+ \bar{g} \\
&=& 2 \alpha  \bar{\g}_V  \y_V^* + 2 (1-\alpha) \g_V^T \e_V 
+ \bar{g} \\
&=& \alpha (OPT_{SDP} - \bar{g}) + 2 (1-\alpha) \bar{\g}_V^T \e_V + \bar{g} \\
&=& \alpha OPT_{SDP} +  (1-\alpha) (2 \bar{\g}_V^T \e_V + \bar{g}).
\end{eqnarray*}
\qed

From a theoretical viewpoint, Lemma~\ref{le:SDP-R} gives the bounds on the expected objective value $E[2 \bar{\g}_V^T \tilde{\y}_V + \bar{g}]$ 
of the randomized algorithm. When we executed preliminary experiments, we observed that the interval between the lower and upper bounds
are not so sharp.
Table~\ref{table:bounds} presents the lower and upper bounds 
and the expected objective value. The dataset we used here is a subset 
of datasets in Section~5.
The first column shows $Z$, the number of genotype candidates. We fix the number of chosen candidates to $N=50$.
% and the threshold for group coancestry to $\theta = 0.015$.
The third and fourth columns are the lower and the upper bounds in  Lemma~\ref{le:SDP-R}, respectively.
The expected objective value is shown in the third column, and it is obtained by generating the random vector $\v$ 
thousand times and taking the average of the thousand trials. 
The fifth column is the optimal value of the SDP relaxation problem~\ref{eq:SDP-R}.

\begin{table}[tbp]
\begin{center}
\caption{Theoretical bounds on the expected values by the randomized algorithm}
\label{table:bounds}
\begin{tabular}{r|r|rrr|r}
\hline
\multicolumn{1}{c|}{Z}  & \multicolumn{1}{c|}{$2\theta$}  & lower bound & expected value & upper bound & $OPT_{SDP}$ \\
\hline
200  & 0.0334 & 16.161 & 25.812 & 30.340 & 25.386 \\
1050 &  0.0627 & 5.075 & 32.305 & 112.600 & 24.938 \\
2045 & 0.0711 & 279.259 & 446.089 & 2007.212 & 438.659 \\
5050 & 0.1081 & 5.775 & 284.965 & 806.205 & 42.786 \\
% 10100 & 0.070094 & 3.04610 & 235.06660 & 1256.11217 & 44.6617 \\
%\hline
%200  & 14.51821 & 22.85170 & 28.07476 & 22.80516 \\
%1050 &  -4.167075 & 27.83562 & 99.85434 & 10.42119 \\
%2045 & 245.8239 & 406.4667 & 1961.100 & 386.1393 \\
%5050 & -5.588337 & 81.12698 & 790.5341 & 15.63910 \\
%5255 & 52.05141 & 645.5305 & 6644.4586 & 322.0503 \\
\hline
\end{tabular}
\end{center}
\end{table}
For the smallest size $Z = 200$, the gap between the lower and the upper bound was not so large.
However, when we tried the larger problems, the gap was getting worse.
In particular, the ratio of the upper bound to the expected objective value 
for the case $Z = 5255$
goes beyond 5.34.

Another aspect in the randomized algorithm is that 
the expected objective value is always larger than $OPT_{SDP}$.
% in contrast to Lemma~\ref{le:relaxation}.
A reason of this unfavorable aspect is that the generated solution $\tilde{\y}_V$ 
is not guaranteed to satisfy the constraint $\y_V^T \A_{VV} \y_V + 2 \bar{\c}_F^T \y_V \le 2 \bar{\theta}$
that corresponds to  $\x^T \A \x \le 2 \theta$ of (\ref{eq:ED}).
Though Theorem~4 of \cite{tseng2003further} estimates the number of randomly generated solutions
required for approximate feasible solutions with high probability,
this cannot be applied to the discussion in this paper,
since the current discussion does not fully satisfy the assumption of the theorem.

Due to this weaker bounds reported in Table~\ref{table:bounds},
we are determined to seek an optimization method 
that can obtain a reasonable solution for practical use.
This motivated us to develop a local search method based on the steep-descent method for discrete convex functions.

%\begin{eqnarray*}
%Var(\bar{\g}_V^T \tilde{\y}_V)
%&=& E[(\sum_{i \in V} g_i \tilde{y}_i)^2] 
%- E[\sum_{i \in V} g_i \tilde{y}_i]^2 \\
%&=& \sum_{i \in V} \sum_{j \in V} g_i g_j 
%E[\tilde{y}_i \tilde{y}_j] 
%- (\sum_{i \in V} g_i  E[\tilde{y}_i])^2 \\
%&=& \sum_{i \in V} \sum_{j \in V} g_i g_j 
%\mbox{arcsin}(Y_{ij}^*)
% - (\sum_{i \in V} g_i  \mbox{arcsin}(Y_{0i}^*))^2
%\end{eqnarray*}
%

\section{Steepest-ascent method}

In contrast to mixed-integer linear programming problems for which
many solvers have been developed,
a principal difficulty in the ED problem (\ref{eq:ED}) arises from 
the nonlinear constraint $\x^T \A \x \le 2 \theta$.
To obtain a sensible solution in a short time,
we embed the violation against this constraint into the objective function as a penalty term
using a penalty weight $\lambda \ge 0$ and focus the following optimization problem
\begin{eqnarray}
\begin{array}{lcl}
\max &:& f_{\lambda}(\x) := \g^T \x - \lambda \max\{\x^T \A \x - 2 \theta, 0\} \\
\mbox{subject to} &:& \x \in \hat{\FC} 
\end{array} \label{eq:penalty}
\end{eqnarray}
where $\hat{\FC} := \left\{ \x \in \Real^Z: \e^T \x = 1, \l \le \x \le \u,  x_1, \ldots, x_Z \in \left\{0, \frac{1}{N}\right\} \right\}$.

We give a validity  of (\ref{eq:penalty}) by  the next lemma
which shows that if we take a large $\lambda$,
this optimization problem with a penalty term (\ref{eq:penalty}) 
is equivalent to the original problem (\ref{eq:ED}) .

\begin{LEMM}\label{le:penalty}
Let $\x(\lambda) \in \Real^Z $ be  an optimal solution of (\ref{eq:penalty}).
There exists a $\hat{\lambda} > 0$ such that  $\x(\lambda)$ is an optimal solution 
of (\ref{eq:ED}) for $\forall \lambda \ge \hat{\lambda}$.
\end{LEMM}

\textbf{Proof:}
Let $\hat{\phi}$ be the optimal value of the following optimization problem;
\begin{eqnarray}
\begin{array}{lclcl}
\hat{\phi} &:=& \min &:& \max\{\x^T \A \x - 2 \theta, 0\} \\
&& \mbox{subject to} &:& \x \in \hat{\FC}. 
\end{array} \label{eq:phi}
\end{eqnarray}
From this definition, $\phi$ can take either zero or a positive number.

If $\hat{\phi} = 0$, the quadratic constraint $\x^T \A \x \le 2 \theta$  holds for $\forall \x \in \hat{\FC}$.
Therefore, this constraint vanishes from (\ref{eq:ED}) and the penalty term % with $\lambda$ 
in (\ref{eq:penalty}) has no effect.
Hence, the two problems (\ref{eq:ED}) and (\ref{eq:penalty}) are equivalent 
% for any $\lambda \in \Real$.
for any $\lambda \ge 0$.

For the case $\hat{\phi} > 0$, 
since $\hat{\FC}$ is composed of a finite number of points, 
the reciprocal number of $\hat{\phi}$ is  a finite number.
Therefore, we can take
 $\hat{\lambda} = \frac{ \max\{ g_i : i = 1, \ldots, Z\} - \min\{ g_i : i = 1, \ldots, Z\} + 1}{\hat{\phi}}$.
To show this by a contradiction, we assume that 
$\x(\lambda)^T \A \x(\lambda) \le 2 \theta$ does not hold  for  $\lambda \ge \hat{\lambda}$. Then, we have 
\begin{eqnarray*}
\g^T \x(\lambda) - \lambda \left( \x(\lambda)^T \A \x(\lambda) - 2\theta \right) 
\le \max\{ g_i : i = 1, \ldots, Z\} - \hat{\lambda} \hat{\phi}
< \min\{ g_i : i = 1, \ldots, Z\}.
\end{eqnarray*}
Here, we used $\e^T \x(\lambda) = 1$ and $\x(\lambda) \ge \0$ since $\x(\lambda) \in \hat{\FC}$.
On the contrary, from the assumption that (\ref{eq:ED}) has a feasible point, the objective value of (\ref{eq:penalty}) at this feasible point 
is at least $\min\{ g_i : i = 1, \ldots, Z\}$. This indicates that $\x(\lambda)$ can not be an optimal solution of 
(\ref{eq:penalty}) if $\x(\lambda)^T \A \x(\lambda) > 2 \theta$.
% , in other words,
%  $\x(\lambda)$ always satisfies the constraint $\x^T \A \x \le 2 \theta$.
Therefore, we can restrict the feasible region of % the penalized problem 
(\ref{eq:penalty}) to the set $\{\x \in \Real^Z : \x^T \A \x \le 2 \theta\} \cap \hat{\FC}$,
and the objective function of (\ref{eq:penalty}) is reduced to $\g^T \x$.
Consequently, the optimal solution of (\ref{eq:penalty}) is also optimal 
for (\ref{eq:ED}).
\qed

Since the computation for $\hat{\phi}$ is almost as hard as the original ED problem,
it is not practical to compute $\hat{\phi}$. In addition, when we maximize 
$f_{\lambda}(\x)$, extremely large $\lambda$ makes the computation numerically unstable.
As an appropriate value for the penalty weight $\lambda$, we make the use of
the Lagrangian multiplier $\lambda_0$ developed in Meuwissen~\cite{MEUWISSEN97};
\begin{eqnarray*}
\lambda_0 := \sqrt{\frac{(\g^T \A^{-1} \g) (\e^T \A^{-1} \e)
              - (\g^T \A^{-1} \e)^2}
             {8\theta (\e^T \A^{-1} \e) - 4}}.
\end{eqnarray*}
This $\lambda_0$ corresponds to the Lagrangian multiplier 
of the constraint $\x^T \A \x = 2 \theta$ in the following optimization problem.
\begin{eqnarray}
\begin{array}{lcl}
\max &:&  \g^T \x  \\
\mbox{subject to} &:& \x^T \A \x = 2\theta,\\
& & \e^T \x = 1.
\end{array}\label{eq:Meuwissen}
\end{eqnarray}
The approach of~\cite{MEUWISSEN97} first solves (\ref{eq:Meuwissen}), 
then applies some heuristic method to obtain a solution of (\ref{eq:ED}).
Therefore, a maximization of $\g^T\x - \lambda_0 (\x^T \A \x - 2 \theta)$
over $\e^T\x = 1$  is a natural derivation when we consider  (\ref{eq:Meuwissen}).
For $f_{\lambda}(\x)$, we often employ $\lambda$ such that $ \lambda \ge \lambda_0$, 
since we put a strong emphasis on the violation with respect to
$\max\{\x^T \A \x - 2 \theta, 0\}$.

We now discuss (\ref{eq:penalty}) from the viewpoint of convex functions.
The function $-f_{\lambda}(\x)$ is a convex function 
in the continuous space $\Real^Z$, since $\A \succeq \O$
and $\lambda \ge 0$. Hence, the problem (\ref{eq:penalty})
can be cast as a minimization of a convex function over a discrete feasible set.

An M-convex function~\cite{murota2003discrete} is a discrete convex function 
defined on a set in which the sum of the elements of 
a feasible point is constant. % From the context of discrete convex functions,
A steepest-descent method  for M-convex functions
was developed in~\cite{murota2004steepest}.
When 
 an M-convex function $f^M (\x)$ with a feasible set $\FC^M$ is given,
the steepest-descent method starts from an initial point $\x^0 \in \FC^M$, 
and finds the next point $\x^1$ from a neighborhood $\NC(\x^0) \subset \FC^M$
that decreases the objective function $f^M (\x)$ with the largest margin. Here, 
$\NC(\x^0) := \{\x + \e_i - \e_j \in \FC^M : i, j = 1, \ldots, n\}$.
In other words, $\x^1$ is chosen so that % it satisfies
$f^M(\x_1) \le f^M(\x)$ for any $\x \in \NC(\x^0)$.
The steepest-descent method continues the search in neighbors,
and it eventually can find a global minimizer 
since any local minimizer is a global minimizer when the objective function $f^M$
is an M-convex function.

Though $-f_{\lambda}(\x)$ is not an M-convex function since it can encompass
multiple local minimizers that are not always global minimizers,
the optimization problem with the penalty term (\ref{eq:penalty}) 
has resemblances to a minimization of an M-convex function.
In particular, the feasible set $\hat{\FC}$ satisfies $\sum_{i=1}^Z x_i = 1$ 
and the function $-f_{\lambda}(\x)$ is a convex function 
in the continuous space $\Real^Z$.
Therefore, we can expect that the steepest-descent method for M-convex functions
will give good direction to solve~(\ref{eq:penalty}).
Furthermore, we can exploit the solution obtained by the conic relaxation problems
in Section~3 to generate a starting point $\x^0$.

When we adjust the steepest-descent method implemented in the software package ODICON~\footnote{\url{http://ist.ksc.kwansei.ac.jp/~tutimura/odicon/index.en.html}}~\cite{TSUCHIMURA} to solve (\ref{eq:penalty}),
we obtain Algorithm~\ref{al:steep}.
Since (\ref{eq:penalty}) is a maximization problem, 
Algorithm~\ref{al:steep} is a steepest-ascent method.

\begin{ALGO}\label{al:steep}
A steep-ascent method with a conic relaxation problem for the optimization problem
with the penalty term arising from the ED problem

\begin{adjustwidth}{1.5em}{}
\begin{enumerate}
\item[Step 1:]  Solve a conic relaxation problem
(\ref{eq:SDP-R}), (\ref{eq:LP-R}) or (\ref{eq:SOCP-R}).
If (\ref{eq:SOCP-R}) is solved, let $\x^*$ be its optimal solution.
For (\ref{eq:LP-R}) and (\ref{eq:SDP-R}), let $\y_V^*$ be its optimal solution
and set $\x^*$ by $\x_V^* := \y_V^*$ and $\x_F^* := \c_F$.
\item[Step 2:] By sorting $\x^*$,
separate $V$ into the two disjoint set $V_{\frac{1}{N}}$  and $V_0$ 
such that $x_i^* \ge x_j^*$ for $i \in V_{\frac{1}{N}}$, $j \in V_0$
and that $|V_{\frac{1}{N}}| = N-p$
(ties are broken arbitrary).
Set the initial point $\x^0 \in \Real^Z$ by
$x_i^0 := \frac{1}{N}$ for $i \in V_{\frac{1}{N}}$,
$x_j^0 := 0$ for $j \in V_0$, and $\x_F^0 := \c_F$.
Set the iteration counter $h:=0$.
\item[Step 3:]
Select the steepest swap $i^h \in V_{\frac{1}{N}}$ and $j^h \in V_0$ such that
\begin{eqnarray*}
f_{\lambda}(\x^h - \frac{1}{N}\e_{i^h} + \frac{1}{N}\e_{j^h}) 
\ge f_{\lambda}(\x^h - \frac{1}{N}\e_i + \frac{1}{N}\e_j) \mbox{ for }
i \in V_{\frac{1}{N}}, j \in V_0.
\end{eqnarray*}
\item [Step 4:]
If there is no improvement, that is $f_{\lambda}(\x^h - \frac{1}{N}\e_{i^h} +\frac{1}{N}\e_{j^h}) \le f_{\lambda}(\x^h)$,
output $\x^h$ as a solution and stop.
\item [Step 5:]
Set $\x^{h+1} := \x^h - \frac{1}{N}\e_{i^h} + \frac{1}{N}\e_{j^h}$.
Swap $i^h$ and $j^h$ by
$V_{\frac{1}{N}} := V_{\frac{1}{N}} \cup \{j^h\} \backslash \{i^h\}$
and $V_0 := V_0 \cup \{i^h\} \backslash \{j^h\}$.
Set $h:= h+1$ and return to Step 3.
\end{enumerate}
\end{adjustwidth}
\end{ALGO}

In Step~2, the number of $\frac{1}{N}$ in $\x^0$ is exactly $N$.
Due to Step~5, this property is kept through the iterations in the algorithm, hence,
the number of $\frac{1}{N}$ in $\x^h$ is also exactly $N$
for any $h \ge 1$.
When no improvement can be found, the algorithm stops by Step~4.

Most computation cost of each iteration in Algorithm~\ref{al:steep}
is consumed at the evaluations of $f_{\lambda}$ in Step~3.
The number of the evaluations is  determined by 
the size of neighbor around $\x^h$, that is, 
$|V_{\frac{1}{N}}| \times |V_0| = (N-p) \times (|V| - (N-p))$.
Therefore, the case $N-p = \frac{|V|}{2}$ requires the heaviest computation cost.
%Furthermore, to reduce the computation cost, we evaluate the difference
%$f_{\lambda}(\x^h - \frac{1}{N}\e_{i^h} +\frac{1}{N}\e_{j^h}) - f_{\lambda}(\x^h)$
%instead $f_{\lambda}(\x^h - \frac{1}{N}\e_{i^h} +\frac{1}{N}\e_{j^h})$ itself.
%This is because
%we need to access all the elements of $\A$ 
%when compute $f_{\lambda}(\x^h - \frac{1}{N}\e_{i^h} +\frac{1}{N}\e_{j^h})$,
%while only the $i^h$th and $j^h$th  rows and columns of $\A$ are involved 
%to evaluate the difference.
Furthermore, to reduce the computation cost, we focus the evaluation of
$(\x^h - \frac{1}{N} \e_{i^h} + \frac{1}{N} \e_{j^h})^T \A 
(\x^h - \frac{1}{N} \e_{i^h} + \frac{1}{N} \e_{j^h})$.
Since ODICON was designed to handle general functions, 
it accessed all the elements of $\A$ for each $i^h$ and $j^h$.
By expanding the part as
$(\x^h)^T \A (\x^h) 
+ \frac{2}{N} \left( - \e_{i^h} + \e_{j^h} \right)^T (\A \x^h)
+ \frac{1}{N^2} \left(A_{i^hi^h} + A_{j^hj^h} - 2 A_{i^hj^h}\right)$,
we evaluate $(\x^h)^T \A (\x^h)$ and $(\A \x^h)$ only once for each iteration 
of Algorithm~\ref{al:steep}.
This saves  95\% of the computation time for Step 3 compared to ODICON.

\section{Numerical results}
In this section, we report numerical results to verify the performance
of the proposed algorithm, Algorithm~\ref{al:steep}.
We implemented Algorithm~\ref{al:steep} with Matlab~R2014b.
We compared the proposed algorithm with GENCONT ~\cite{MEUWISSEN97},
the branch-and-bound method implemented in OPSEL~2.0~\cite{mullin2014opsel, mullin2016using}, and IBM CPLEX~12.62.
We used an Windows PC with Core i7 3770K (3.5 GHz) and 32 GB memory space for cases.
Only when the 32 GB memory space was not enough, we used a Linux server 
with Opteron 4386 (3.10 GHz) and 128 GB memory space.
To solve the LP problem~(\ref{eq:LP-R}), the SOCP problem~(\ref{eq:SOCP-R}), and the SDP problem~(\ref{eq:SDP-R}),
we employed CPLEX, ECOS~\cite{DOMAHIDI13}, and SDPT-3~\cite{TODD99}, respectively.
% The stopping tolerances in SDPT-3 were relaxed to $10^{-4}$ from $10^{-8}$.
For the steepest-ascent method, we set $\lambda = 2\lambda_0$ as the penalty weight in
the function $f_{\lambda}(\x)$ of ~(\ref{eq:penalty}).

The data tested in the numerical experiments of this paper 
are % same one used in Section~3;
practical datasets of pine orchards available at
the Dryad Digital Repository\footnote{\url{http://dx.doi.org/10.5061/dryad.9pn5m}}
and datasets generated by a simulation software package~\cite{mullin2010using}.

Tables~\ref{table:result50} and \ref{table:result100} present the comparison of 
the three conic relaxation approaches.
% GENCONT, OPSEL,
%  and the proposed algorithm~(Algorithm~\ref{al:steep}).
The number $N$ of chosen genotype is set 50 and 150 in
Table~\ref{table:result50} and Table~\ref{table:result100}, respectively.
% and the threshold is set $2 \theta = 0.0300$ 
% in Tables~\ref{table:result50}, \ref{table:result100}, \ref{table:result150}, respectively.
The first column in the tables shows the name of algorithms; 
for example, CR (LP) is the result of conic relaxation problem (in this case, an LP problem) and 
SA (LP) is the result after the application of Algorithm~\ref{al:steep} starting from the solution of CR (LP).
The names for SOCP and SDP are indicated with the same rule.
%  and these names are listed in Table~\ref{table:name}.
For CR (LP)-s and SA (LP)-s, we applied Remark~\ref{re:LP} to (\ref{eq:LP-R})
and obtain its solution by sorting $\g_V$.
The second column is $Z$, the number  of genotype candidates, 
while the third column is $2 \theta$.
The fourth and fifth columns are the objective value $\g^T \x$ and the value $\x^T \A \x$
of each algorithm. For the CR rows, these two values were evaluated at $\x^*$,
the solution at Step~1 of Algorithm~\ref{al:steep}, and for the SA rows,
they were computed with the output solution $\x^h$ at Step~4.
%When an algorithm output $\y_V$ as a solution (for example, SA (SDP)),
%we revert it to $\x \in \Real^Z$ by setting $\x_V := \frac{\y_V + \e_V}{2}$ and $\x_F := \c_F$.
The sixth column is the iteration number of Algorithm~\ref{al:steep}.
The seventh column is the value of $f_{\lambda}(\x^h)$,
where $h$ is the iteration number indicated in the sixth column.
For the CR row, note that $\x^*$ must satisfy the quadratic 
constraint $\x^T \A \x \le 2 \theta$,
but $\x^0$ does not always satisfy it.
In contrast, for the SA rows, $f_{\lambda}(\x^h)$ is given at the final  solution of Step~4.
The last column is the computation time in seconds.
Since SA (LP) uses the result of CR (LP), the computation time of SA (LP) is the sum of
the computation time of CR (LP) and the steep-ascent method.
In a similar way, SA (SOCP) and SA (SDP) are also the sums.

\begin{table}[tbp]
\begin{center}
\caption{The comparison of the convex relaxation approaches ($N = 50$)}\label{table:result50}
{\small 
\
\begin{tabular}{l|c|c|c|c|r|r|r}
\hline
Algorithm & $Z$  & $2 \theta$ & $\g^T \x$ & $\x^T \A \x$ & iter & \multicolumn{1}{|c|}{$f_{\lambda}(\x)$} & time (s) \\
\hline
% @@@@@@@@@@@@@ Final Result TeX @@@@@@@@@
%---  Z = 200, N = 50, 2*theta = 0.033430, 2*theta_relax = 0.031600
CR (LP)  & \multirow{8}{*}{200} & \multirow{8}{*}{0.0334} & 28.068 &  0.0574 & 0 & -133.895  & 0.07 \\
SA (LP)    &  & & 25.029 &  0.0334  & 21 & 25.029  & 0.10 \\
CR (LP)-s  &  & & 28.068 &  0.0574 & 0 & -133.895  & 0.01 \\
SA (LP)-s  & & & 25.029 &  0.0334 & 21 & 25.029  & 0.04 \\
CR (SOCP)  & & & 26.156 &  0.0334 & 0 &  -41.484  & 0.02 \\
SA (SOCP)  & & & 25.090 &  0.0334  & 13 & 25.090  & 0.06 \\
CR (SDP)  & & & 25.386 &  0.0321 & 0 & 18.978  & 1.29 \\
SA (SDP)  & & & 25.207 &  0.0334 & 4 & 25.207  & 1.30 \\
% SDP_lower = 16.16101, SDP_ave = 25.81230, SDP_upper = 30.34043
% @@@@@@@@@@@@@ Final Result TeX @@@@@@@@@
\hline
% @@@@@@@@@@@@@ Final Result TeX @@@@@@@@@
%---  Z = 1050, N = 50, 2*theta = 0.062748, 2*theta_relax = 0.030000
CR (LP)  & \multirow{8}{*}{1050} & \multirow{8}{*}{0.0627} & 30.754 &  0.1362 & 0 & -198.235  & 0.37 \\
SA (LP)    &  & & 22.707 &  0.0627  & 23 & 22.707  & 0.51 \\
CR (LP)-s  &  & & 30.754 &  0.1362 & 0 & -198.235  & 0.01 \\
SA (LP)-s  & & & 22.707 &  0.0627 & 23 & 22.707 & 0.15 \\
CR (SOCP)  & & & 25.284 &  0.0627 & 0 & 19.621  & 0.08 \\
SA (SOCP)  & & & 24.831 &  0.0627  & 2 & 24.831  & 0.09 \\
CR (SDP)  & & & 24.938 &  0.0617 & 0 & 24.721  & 27.94 \\
SA (SDP)  & & & 24.846 &  0.0627 & 2 & 24.846 & 27.96 \\
% SDP_lower = 5.07464, SDP_ave = 32.30552, SDP_upper = 112.60014
% @@@@@@@@@@@@@ Final Result TeX @@@@@@@@@
\hline
% @@@@@@@@@@@@@ Final Result TeX @@@@@@@@@
%---  Z = 2045, N = 50, 2*theta = 0.071083, 2*theta_relax = 0.029400
CR (LP)  & \multirow{8}{*}{2045} & \multirow{8}{*}{0.0711} & 504.217 &  0.4566 & 0 & -26197.137  & 1.16 \\
SA (LP)    &  & & 414.591 &  0.0710  & 32 & 414.591  & 1.47 \\
CR (LP)-s  &  & & 504.217 &  0.4566 & 0 & -26197.137  & 0.01 \\
SA (LP)-s  & & & 414.591 &  0.0710 & 32 & 414.591  & 0.32 \\
CR (SOCP)  & & & 439.353 &  0.0711 & 0 & 293.122  & 0.06 \\
SA (SOCP)  & & & 438.386 &  0.0710  & 2 & 438.386  & 0.09 \\
CR (SDP)  & & & 438.659 &  0.0706 & 0 & 438.457  & 145.57 \\
SA (SDP)  & & & 438.457 &  0.0710 & 1 & 438.457  & 145.59 \\
% SDP_lower = 279.25875, SDP_ave = 446.08925, SDP_upper = 2007.21232
% @@@@@@@@@@@@@ Final Result TeX @@@@@@@@@
\hline
%---  Z = 5050, N = 50, 2*theta = 0.108061, 2*theta_relax = 0.031625
CR (LP)  & \multirow{8}{*}{5050} & \multirow{8}{*}{0.1081} & 57.630 &  0.3672 & 0 & -1185.866  & 10.17 \\
SA (LP)    &  & & 38.696 &  0.1080  & 23 & 38.696  & 11.17 \\
CR (LP)-s  &  & & 57.630 &  0.3672 & 0 & -1185.866  & 0.01 \\
SA (LP)-s  & & & 38.696 &  0.1080 & 23 & 38.696  & 0.98 \\
CR (SOCP)  & & & 43.036 &  0.1081 & 0 & 42.456  & 0.21 \\
SA (SOCP)  & & & 42.691 &  0.1080  & 3 & 42.691  & 0.37 \\
CR (SDP)  & & & 42.786 &  0.0980 & 0 & 41.327  & 2221.22 \\
SA (SDP)  & & & 42.431 &  0.1080 & 3 & 42.431  & 2221.40 \\
% SDP_lower = 5.77446, SDP_ave = 284.96506, SDP_upper = 806.20523
% @@@@@@@@@@@@@ Final Result TeX @@@@@@@@@
\hline
% @@@@@@@@@@@@@ Final Result TeX @@@@@@@@@
%---  Z = 10100, N = 50, 2*theta = 0.070094, 2*theta_relax = 0.000000
CR (LP)  & \multirow{8}{*}{10100} & \multirow{8}{*}{0.0701} & 62.377 &  0.2368 & 0 & -1305.4682  & 46.84 \\
SA (LP)    &  & & 41.284 &  0.0701  & 32 & 41.284 & 49.87 \\
CR (LP)-s  &  & & 62.377 &  0.2368 & 0 & -1305.468  & 0.01 \\
SA (LP)-s  & & & 41.284 &  0.0701 & 32 & 41.284  & 3.29 \\
CR (SOCP)  & & & 47.445 &  0.0701 & 0 & 21.094  & 0.54 \\
SA (SOCP)  & & & 46.568 &  0.0701  & 2 & 46.568  & 0.87 \\
CR (SDP)  & & & 21.265 &  0.0545 & 0 & 13.369  & 5577.80$\dagger$ \\
SA (SDP)  & & & 44.662 &  0.0701 & 45 & 44.662 & 5582.46$\dagger$ \\
% SDP_lower = 3.04610, SDP_ave = 235.06660, SDP_upper = 1256.11217
% @@@@@@@@@@@@@ Final Result TeX @@@@@@@@@
\hline
% @@@@@@@@@@@@@ Final Result TeX @@@@@@@@@
%---  Z = 15222, N = 50, 2*theta = 0.038808, 2*theta_relax = 0.000000
CR (LP)  & \multirow{8}{*}{15222} & \multirow{8}{*}{0.0388} & 603.783 &  0.4568 & 0 & -67047.589  & 129.55 \\
SA (LP)    &  &&  438.791 &  0.0388  & 42 & 438.791  & 139.03 \\
CR (LP)-s  &  & & 603.783 &  0.4568 & 0 & -67047.589  & 0.01 \\
SA (LP)-s  & & & 438.791 &  0.0388 & 42  & 438.791 & 6.45 \\
CR (SOCP)  & & & 468.367 &  0.0388 & 0 & -1042.485 & 0.99 \\
SA (SOCP)  & & & 460.769 &  0.0388  & 9 & 460.769  & 2.56 \\
CR (SDP)  & & & 288.739 &  0.0195 & 0 & 314.493  & 17433.38$\dagger$ \\
SA (SDP) & & & 460.409 &  0.0388 & 43 & 460.409 & 17441.93$\dagger$ \\
%---  Z = 15222, N = 50, 2*theta = 0.038808, 2*theta_relax = 0.000000
% @@@@@@@@@@@@@ Final Result TeX @@@@@@@@@
\hline
\end{tabular}
} % end of \small
\end{center}
$\dagger$ indicates numerical instability 
\end{table}

\begin{table}[tbp]
\begin{center}
\caption{The comparison of the convex relaxation approaches ($N = 100$)}\label{table:result100}
{\small 
\begin{tabular}{l|c|c|c|c|r|r|r}
\hline
Algorithm & $Z$  & $2 \theta$ & $\g^T \x$ & $\x^T \A \x$ & iter & \multicolumn{1}{|c|}{$f_{\lambda}(\x)$}  &  time (s) \\
\hline
% @@@@@@@@@@@@@ Final Result TeX @@@@@@@@@
%---  Z = 200, N = 100, 2*theta = 0.025817, 2*theta_relax = 0.031600
CR (LP)  & \multirow{8}{*}{200} & \multirow{8}{*}{0.0258} & 24.654 &  0.0304 & 0 & -22.392 & 0.01 \\
SA (LP)    &  & & 23.355 &  0.0258  & 18 & 23.355 & 0.07 \\
CR (LP)-s  &  & & 24.654 &  0.0304 & 0 & -22.392 &  0.01 \\
SA (LP)-s  & & & 23.355 &  0.0258 & 18 & 23.355 &  0.06 \\
CR (SOCP)  & & & 24.015 &  0.0258 & 0 & 0.950 &  0.01 \\
SA (SOCP)  & & & 23.412 &  0.0258  & 13 & 23.412 & 0.08 \\
CR (SDP)  & & & 23.640 &  0.0255 & 0 & 21.783 &  0.82 \\
SA (SDP)  & & & 23.521 &  0.0258 & 5 & 23.521 &  0.84 \\
% SDP_lower = 15.04972, SDP_ave = 23.56180, SDP_upper = 24.78188
% @@@@@@@@@@@@@ Final Result TeX @@@@@@@@@
\hline
% @@@@@@@@@@@@@ Final Result TeX @@@@@@@@@
%---  Z = 1050, N = 100, 2*theta = 0.053925, 2*theta_relax = 0.030000
CR (LP)  & \multirow{8}{*}{1050} & \multirow{8}{*}{0.0539} & 27.637 &  0.1214 & 0 & -208.680  & 0.39 \\
SA (LP)    &  & & 19.805 &  0.0539  & 42 & 19.805  & 0.98 \\
CR (LP)-s  &  & & 27.637 &  0.1214 & 0 & -208.680  & 0.01 \\
SA (LP)-s  & & & 19.805 &  0.0539 & 42 & 19.805 & 0.590 \\
CR (SOCP)  & & & 22.432 &  0.0539 & 0 & 18.537  & 0.08 \\
SA (SOCP)  & & & 22.321 &  0.0539  & 5 & 22.321  & 0.15 \\
CR (SDP)  & & & 22.358 &  0.0537 & 0 & 22.242  & 30.49 \\
SA (SDP)  & & & 22.324 &  0.0539 & 3 & 22.324  & 30.54 \\
% SDP_lower = 3.43185, SDP_ave = 23.93306, SDP_upper = 63.16898
% @@@@@@@@@@@@@ Final Result TeX @@@@@@@@@
\hline
% @@@@@@@@@@@@@ Final Result TeX @@@@@@@@@
%---  Z = 2045, N = 100, 2*theta = 0.062824, 2*theta_relax = 0.029400
CR (LP)  & \multirow{8}{*}{2045} & \multirow{8}{*}{0.0628} & 478.114 &  0.4219 & 0 & -26349.725  & 1.17 \\
SA (LP)    &  & & 406.348 &  0.0628  & 65 & 406.348  & 5.00 \\
CR (LP)-s  &  & & 478.114 &  0.4219 & 0 & -26349.725  & 0.01 \\
SA (LP)-s  & & & 406.348 &  0.0628 & 65 & 406.348  & 3.78 \\
CR (SOCP)  & & & 421.696 &  0.0628 & 0 & 197.423  & 0.07 \\
SA (SOCP)  & & & 421.113 &  0.0627  & 3 & 421.113  & 0.25 \\
CR (SDP)  & & & 421.497 &  0.0627 & 0 & 364.014  & 165.33 \\
SA (SDP)  & & & 421.425 &  0.0628 & 2 & 421.425  & 165.46 \\
% SDP_lower = 268.33355, SDP_ave = 426.98624, SDP_upper = 1181.10971
% @@@@@@@@@@@@@ Final Result TeX @@@@@@@@@
\hline
%---  Z = 5050, N = 100, 2*theta = 0.099463, 2*theta_relax = 0.031625
CR (LP)  & \multirow{8}{*}{5050} & \multirow{8}{*}{0.0994} & 54.903 &  0.3355 & 0 & -1137.701  & 10.35 \\
SA (LP)    &  & & 36.509 &  0.0994  & 50 & 36.509 & 17.44 \\
CR (LP)-s  &  & & 54.903 &  0.3355 & 0 & -1137.701  & 0.01 \\
SA (LP)-s  & & & 36.509 &  0.0994 & 50 & 36.509  & 7.03 \\
CR (SOCP)  & & & 40.769 &  0.0995 & 0 & 15.408 & 0.25 \\
SA (SOCP)  & & & 40.629 &  0.0995  & 3 & 40.629  & 0.71 \\
CR (SDP)  & & & 40.711 &  0.0992 & 0 & 31.447  & 2164.49 \\
SA (SDP)  & & & 40.690 &  0.0994 & 2 & 40.690  & 2164.85 \\
% SDP_lower = 10.37304, SDP_ave = 56.74496, SDP_upper = 421.53648
% @@@@@@@@@@@@@ Final Result TeX @@@@@@@@@
\hline
% @@@@@@@@@@@@@ Final Result TeX @@@@@@@@@
%---  Z = 10100, N = 100, 2*theta = 0.060999, 2*theta_relax = 0.000000
CR (LP)  & \multirow{8}{*}{10100} & \multirow{8}{*}{0.0610} & 60.347 &  0.2245 & 0 & -1395.786  & 46.71 \\
SA (LP)    &  & & 39.911 &  0.0610  & 62 & 39.911 & 68.54 \\
CR (LP)-s  &  & & 60.347 &  0.2245 & 0 & -1395.786  & 0.01 \\
SA (LP)-s  & & & 39.911 &  0.0610 & 62 & 39.911  & 21.51 \\
CR (SOCP)  & & & 44.819 &  0.0610 & 0 & 10.608  & 0.70 \\
SA (SOCP)  & & & 44.522 &  0.0610  & 7 & 44.522  & 3.12 \\
CR (SDP)  & & & 21.374 &  0.0532 & 0 & 18.463  & 6750.06$\dagger$ \\
SA (SDP)  & & & 42.810 &  0.0610 & 66 & 42.810  & 6773.05$\dagger$ \\
% SDP_lower = 3.11489, SDP_ave = 138.19201, SDP_upper = 635.72525
% @@@@@@@@@@@@@ Final Result TeX @@@@@@@@@
%\hline
\hline
% @@@@@@@@@@@@@ Final Result TeX @@@@@@@@@
%---  Z = 15222, N = 100, 2*theta = 0.030044, 2*theta_relax = 0.000000
CR (LP)  & \multirow{8}{*}{15222} &  \multirow{8}{*}{0.0300} &575.227 &  0.4318 & 0 & -74482.507  & 128.21 \\
SA (LP)    &  & & 408.725 &  0.0300  & 90 & 408.725 & 185.33 \\
CR (LP)-s  &  & & 575.227 &  0.4318 & 0 & -74482.507  & 0.02 \\
SA (LP)-s  & & & 408.725 &  0.0300 & 90 & 408.725 & 49.14 \\
CR (SOCP)  & & & 444.730 &  0.0300 & 0 & -11.0395 & 1.05 \\
SA (SOCP)  & & & 441.438 &  0.0300  & 6 & 441.438 & 4.72 \\
CR (SDP)  & & & 309.173 &  0.0228 & 0 & -65291.694  & 19467.30$\dagger$ \\
SA (SDP)  & & & 406.266 &  0.0300 & 92 & 406.266  & 19525.52$\dagger$ \\
%---  Z = 15222, N = 50, 2*theta = 0.038808, 2*theta_relax = 0.000000
% @@@@@@@@@@@@@ Final Result TeX @@@@@@@@@
\hline
\end{tabular}
} % end of \small
\end{center}
$\dagger$ indicates numerical instability 
\end{table}

For the smallest case 
$Z = 200$ and $N = 50$, the SDP relaxation problem attains a remarkable result.
Since the solution in SA (SDP) satisfies the constraint $\x^T \A \x \le 2 \theta$,
this solution is a feasible solution of the original ED problem~(\ref{eq:ED}).
In addition, it holds that $OPT_{ED} \le OPT_{SDP}$ from Lemma~\ref{le:relaxation}.
Therefore, we know $25.207 \le OPT_{ED} \le 25.386$
and we obtain the optimal value of the ED problem 
up to an error $0.710\%$.
Since this error is much better than the theoretical bounds discussed in Lemma~\ref{le:SDP-R},
the combination of the SDP relaxation and the steepest-ascent method performs
very well in this case.

In addition, we can make sure that the objective values of CR (LP) 
and CR (LP)-s are same, and 
this indicates that Assumption~\ref{as:A} holds in the numerical tests. 
As a result, we can obtain the solution of 
the LP relaxation problem~(\ref{eq:LP-R}) without solving it as an LP problem,
as noted in Remark~\ref{re:LP}.

From Tables~\ref{table:result50} and \ref{table:result100}, 
we observe in the cases $Z \le 5050$ that 
$OPT_{SDP} \le OPT_{SOCP} \le OPT_{LP}$  from the values of  $\g^T\x$
in the CR~(LP), CR~(SOCP), and CR~(SDP) rows.
This result supports the validity of Lemma~\ref{le:relaxation}.
However, we also observe for large instances $Z \ge 10100$ 
that $\g^T \x$ of CR (SDP) is lower than that of CR (SOCP).
% This phenomenon is inconsistent with Lemma~\ref{le:relaxation},
A principal reason of this inconsistent phenomenon is a premature termination of SDPT-3.
These inaccurate values of CR (SDP) were mainly caused by
the lack of interior-feasible points in~(\ref{eq:SDP-R}); see Remark~\ref{re:SDP}.
%To overcome this difficulty, we also tried 
%a semi-smooth Newton-CG augmented Lagrangian method
%implemented  SDPNAL+~\cite{Yang2015}.
%However, its convergence rate was slow for (\ref{eq:SDP-R}), 
%and it could not attain enough numerical accuracy.
The SOCP relaxation problem~(\ref{eq:SOCP-R}) provides
highly numerical stability compared to the SDP relaxation 
problem~(\ref{eq:SDP-R}). This can be regarded as an advantage 
of the SOCP relaxation approach. 
% extremely large or extremely small elements in $\Y_{VV}$, as pointed out in
% \cite{takeda2000towards}.

As a next viewpoints, % We also observe in these tables that 
the violations of the solution generated 
by the steep-ascent method against $\x^T \A \x \le 2 \theta$ are remarkably small.
This is mainly because we set $\lambda$ large enough 
based on the Lagrangian multiplier $\lambda_0$.
Therefore, the maximization of the function with the penalty term (\ref{eq:penalty}) 
can provide a suitable solution for the ED problem~(\ref{eq:ED}).

% However, this is a rare case. 
% A more difficult situation can be found in  the case of 
%$Z = 1050$, where the discrepancy between $\g^T \x$ of CR~(SDP) and SA~(SDP) 
%is $6.67\%$. Furthermore, $\x^T \A \x$ in SA~(SDP) is slightly larger than $2 \theta$.
%This result accords with Lemma~\ref{le:SDP-R}, 
%which implied that  the theoretical bounds of the SDP relaxation are not always sharp.
%The difficulty of this problem can be also detected from the result of GENCONT in Table (We make this table later);
%the number of chosen genotypes there by GENCONT is 43, not $N = 50$. 

A comparison of  the results of the steepest-ascent methods that starts
from the three conic relaxation indicates that 
if SDPT-3 obtained sensible solutions ($Z \le 5050$),
the output objective values $\g^T \x$ of SA~(SDP) and SA~(SOCP) 
were close to each other, but much higher than that of SA~(LP).
This implies that the SOCP relaxation problem
and the SDP relaxation problem provided good starting points
for the steepest-ascent method.
This point is also indicated in the iteration number of the steepest-ascent methods.
The iteration numbers of SA~(SDP) and SA~(SOCP) 
are much less than that of SA~(LP). 
For example, when $Z = 2045$ and $N = 100$,
SA~(SDP) and SA~(SOCP) required only two and three iterations, respectively,
while SA~(LP) required 65 iterations.
Therefore, we can infer that the solutions of 
the SDP relaxation problem
and the SOCP relaxation problem are close to local maximizer of $f_{\lambda}(\x)$.

When we move our focus from the solution quality to the computation time,
the computation of SA~(SOCP) is much shorter than SA~(SDP).
This was mainly because we can aggressively exploit the structure  of 
the Wright numerator matrix $\A$
through the SOCP approach discussed in \cite{yamashita2015efficient}.
Since the objective values $\g^T \x$ of SA~(SOCP) and SA~(SDP) are competitive,
SA~(SOCP) can be considered as the most efficient approach
among the three SA~(LP), SA~(SOCP), and SA~(SDP).

Judging from these results, we chose the SOCP relaxation problem for generating the starting point 
when we compare Algorithm~\ref{al:steep} 
against the existing approaches, GENCONT, OPSEL, and CPLEX.
Since OPSEL and CPLEX utilized the brand-and-bound framework, 
we can expect their solution are close to the real optimal solution of 
the ED problem~(\ref{eq:ED}).

Tables~\ref{table:compare50} and \ref{table:compare100} present the comparison of 
Algorithm~\ref{al:steep}, GENCONT and OPSEL (Version 2.0), CPLEX (Version 12.62)
for the numbers of selected genotypes $N = 50$ and $N = 100$.
We tried to execute GENCONT2, a new version of GENCONT,
but its binary file did not work on our computational environment. Therefore, we used GENCONT1 for the comparison.
% In these tables, we also add the results of  CR~(SOCP)-r and SA~(SOCP)-r.  
% Since the solution of GENCONT often violates the constraint $\x^T \A \x \le 2 \theta$
% slightly, CR~(SOCP)-r and SA~(SOCP)-r set $2 \theta$ to the value $\x^T \A \x$ 
% obtained with the solution of GENCONT
% so that we can simply compare the results of GENCONT and the proposed method
% by focusing their $\g^T \x$.
In the tables, we evaluate $f_{\lambda}(\x)$ in seventh column at the solution obtained 
from each algorithm
and the eighth column reports 
the number of chosen genotypes $|\{ i : x_i > 0\}|$.
For OPSEL and CPLEX, we set the time limit as three hours and 
the tolerance gap as $1\%$.
When OPSEL and CPLEX reached the time limit,
we obtained a feasible solution from OPSEL, 
but we could not extract sensible solutions from CPLEX.
GENCONT could not solve large instances ($Z = 10100$ and $Z = 15222$)
due to out of memory.

\begin{table}[tbp]
\begin{center}
\caption{The comparison of GENCONT, OPSEL, CPLEX, and the proposed algorithm ($N = 50$)}\label{table:compare50}
{\small 
\begin{tabular}{l|c|c|c|c|c|c|r}
%\hline
%\multicolumn{7}{c}{$N = 50$} \\
\hline
Algorithm & $Z$  & $2{\theta}$  & $ \g^T \x$ & $\x^T \A \x$ &  $f_{\lambda}(\x)$  & \#chosen & time (s) \\
\hline
GENCONT & \multirow{4}{*}{200} & \multirow{4}{*}{0.0334} & 25.290 & 0.0342 & 20.087  & 50 & 0.06 \\
OPSEL & & & 25.191 & 0.0334 & 25.191 & 50 & 1779.13 \\
CPLEX  &  & & 25.190 & 0.0334 & 25.190 & 50 & 4270.77 \\
SA (SOCP)  & &  & 25.090 & 0.0334 & 25.090 & 50 & 0.06 \\
% \cline{3-3}
% SA(SOCP)-r  & & 0.0342 & 25.3720 &  0.0342 &  25.3720 & 50 & 0.038 \\
\hline
GENCONT & \multirow{4}{*}{1050} & \multirow{4}{*}{0.0627} & 24.983 &  0.0627 & 24.983  & 48 & 7.91 \\
OPSEL && &  24.858 & 0.0627 & 24.858 & 50 & $>$ 10800  \\
CPLEX  & & & \multicolumn{4}{c|}{Cannot obtain a sensible solution in 3 hours} & $>$ 10800 \\
SA (SOCP)  & &  & 24.831 & 0.0627 & 24.831 & 50 & 0.09 \\
% \cline{3-3}
% SA(SOCP)-r  & & 0.062748 & 24.8314 &  0.0627  & 24.8314 & 50 & 0.093 \\
\hline
GENCONT & \multirow{4}{*}{2045} & \multirow{4}{*}{0.0711} & 437.049 &  0.0694 & 437.049  & 50 & 88.46 \\
OPSEL & & & 435.826 & 0.0692 & 435.826 & 50 & 16.08 \\
CPLEX  & & & 436.213 & 0.0680 &436.212 & 50 & 0.37 \\
SA (SOCP)  & &  & 438.386 & 0.0710 & 438.386 & 50 & 0.09 \\
% \cline{3-3}
% SA(SOCP)-r  & &  0.0694 & 437.1920 &  0.0692 &  437.1920 & 50 & 0.098\\
\hline
GENCONT & \multirow{4}{*}{5050} & \multirow{4}{*}{0.1081} & 42.780 &  0.1089 & -306.701  & 50 & 1769.72 \\
OPSEL & & & 42.702 & 0.1081 & 42.702 & 50 & $>$ 10800 \\
CPLEX  & & & 42.456 & 0.1066 & 42.456 & 50 & 2.02 \\
SA (SOCP)  & &  & 42.691 & 0.1080 & 42.691 & 50 & 0.37 \\
% \cline{3-3}
% SA(SOCP)-r  & & 0.1089 & 42.8211 &  0.1088 &  42.8211 & 50 & 0.350 \\
\hline
GENCONT & \multirow{4}{*}{10100} & \multirow{4}{*}{0.0701} & \multicolumn{5}{c}{Out of memory} \\
OPSEL & & & 46.252 & 0.0700 & 46.252 & 50 & $>$ 10800 \\
CPLEX  & & & \multicolumn{4}{c|}{Cannot obtain a sensible solution in 3 hours} & $>$ 10800 \\
SA (SOCP)  & &  & 46.568 & 0.0701 & 46.568 & 50 & 0.87 \\
\hline
GENCONT & \multirow{4}{*}{15222} & \multirow{4}{*}{0.0388} & \multicolumn{5}{c}{Out of memory} \\
OPSEL & & & 459.040 & 0.0388 & 459.040 & 50 & $>$ 10800 \\
CPLEX  & & & 459.135 & 0.0386 & 459.135 & 50 & 39.20 \\
SA (SOCP)  & & & 460.769 & 0.0388 & 460.769 & 50 & 2.56 \\
\hline
\end{tabular}
} % end of \small
\end{center}
\end{table}

\begin{table}[tbp]
\begin{center}
\caption{The comparison of GENCONT, OPSEL, CPLEX, and the proposed algorithm ($N = 100$)}\label{table:compare100}
{\small 
\begin{tabular}{l|c|c|c|c|c|c|r}
%\hline
%\multicolumn{7}{c}{$N = 50$} \\
\hline
Algorithm & $Z$  & $2{\theta}$  & $ \g^T \x$ & $\x^T \A \x$ &  $f_{\lambda}(\x)$  & \#chosen & time (s) \\
\hline
GENCONT & \multirow{4}{*}{200} & \multirow{4}{*}{0.0258} & 23.640 &  0.0261 & 21.253  & 100 & 0.07 \\
OPSEL & & & 23.551 & 0.0258 & 23.551 & 100 & $>$ 10800 \\
CPLEX  &  & & 23.508 & 0.0258 & 23.508& 100 & 1.42  \\
SA (SOCP)  & &  & 23.412 & 0.0258 & 23.412 & 100 & 0.08 \\
% \cline{3-3}
% SA(SOCP)-r  & &  0.0261 & 23.6174 &  0.0260 &  23.6174 & 100 & 0.041 \\
\hline
GENCONT & \multirow{4}{*}{1050} & \multirow{4}{*}{0.0539} & 22.749 &  0.0539 &  22.749  & 91 & 9.63 \\
OPSEL & & & 22.275 & 0.0539 & 22.275 & 100 & 304.89 \\
CPLEX  & & & \multicolumn{4}{c|}{Cannot obtain a sensible solution in 3 hours} & $>$ 10800 \\
SA (SOCP)  & &  & 22.321 &  0.0539  & 22.321 & 100 & 0.15 \\
% \cline{3-3}
% SA(SOCP)-r  & &  0.0539 & 22.3211 &  0.0539 &  22.3211 & 100 & 0.146  \\
\hline
GENCONT & \multirow{4}{*}{2045} & \multirow{4}{*}{0.0628} & 421.005 &  0.0632 & 392.893   & 100 & 105.40 \\
OPSEL & & & 419.600 & 0.0613 & 419.600 & 100 & 6.85 \\
CPLEX  & & & 420.748 & 0.0619 & 420.748 & 100 & 0.41 \\
SA (SOCP)  & &  & 421.113 &  0.0627  & 421.113 & 100 & 0.25 \\
% \cline{3-3}
% SA(SOCP)-r  & &  0.0632 & 421.7364 &  0.0632 & 421.7364 & 100 & 0.188 \\
\hline
GENCONT & \multirow{4}{*}{5050} & \multirow{4}{*}{0.0995} & 40.692 &  0.0997 & 40.692  & 100 & 1940.43  \\
OPSEL & & & 40.468 & 0.0994 & 40.468 & 100 & 50.46 \\
CPLEX  & & & \multicolumn{4}{c|}{Cannot obtain a sensible solution in 3 hours} & $>$ 10800  \\
SA (SOCP)  & &  & 40.629 &  0.0995  & 40.629 & 100 & 0.71 \\
% \cline{3-3}
% SA(SOCP)-r  & & 0.0997 & 40.6694 &  0.0997 &  40.6694 & 100 & 0.841  \\ 
\hline
GENCONT & \multirow{4}{*}{10100} & \multirow{4}{*}{0.0610} & \multicolumn{5}{c}{Out of memory} \\
OPSEL & & & 44.467 & 0.0696 & 44.467 & 100 & $>$ 10800 \\
CPLEX  & & & \multicolumn{4}{c|}{Cannot obtain a sensible solution in 3 hours} & $>$ 10800 \\
SA (SOCP)  & &  & 44.522 &  0.0610  & 44.522 & 100 & 3.12 \\
\hline
GENCONT & \multirow{4}{*}{15222} & \multirow{4}{*}{0.0300} & \multicolumn{5}{c}{Out of memory} \\
OPSEL & & & 441.770 & 0.0300 & 441.770 & 100 & $>$ 10800\\
CPLEX  & & & 440.996 & 0.0290 & 440.99640 & 100 & 7.45 \\
SA (SOCP)  & & & 441.438 &  0.0300  & 441.438 & 100 & 4.72 \\
\hline
\end{tabular}
} % end of \small
\end{center}
\end{table}

From the numerical results in Tables~\ref{table:compare50} and \ref{table:compare100}, 
the objective values of SA~(SOCP) are close to those of GENCONT, OPSEL,
and CPLEX. Since the cost vector $\g$ in the objective function is usually generated 
from a statistical procedure, the discrepancy in the objective values 
make a little difference for practical use.
However, GENCONT sometimes failed to satisfy the constraints;
the quadratic constraint $\x^T \A \x \le 2 \theta$ was violated 
in the $Z = 200$ or $Z = 5050$ cases, and 
the number of the chosen genotypes did not match the input $N$.
Therefore, the quality of SA (SOCP) was superior to that of GENCONT.

% and SA~(SOCP) attained this values with much less computation time than OPSEL.
% In addition, SA~(SOCP)-r attained higher  $\g^T \x$ than GENCONT.
% This means that if we admit the same violation to Algorithm~\ref{al:steep} as GENCONT,

From the viewpoint of  the computation time, of SA~(SOCP) is much faster than GENCONT.
In particular, for the case $Z = 5050$, 
 SA~(SOCP) used less than one seconds, while GENCONT required 1700 seconds.
%Another viewpoint is that GENCONT failed to select exactly 50 genotypes for the case $Z = 1050$,
%while the other methods selected 50 genotypes successfully.
The computation times of the branch-and-bound framework were 
unpredictable. In $N = 50$, OPSEL and CPLEX required longer computation 
for a small problem $Z = 200$ than for a large problem $Z = 2045$.
It is very difficult to estimate the computation time required by
OPSEL and CPLEX in advance due to a nature of the branch-and-bound framework.
In contrast, SA~(SOCP) consumed longer computation time for larger problems
and this property is favorable for practical use.

\section{Conclusion and Future Directions}
In this paper, we introduced the conic relaxation approach based on
LP, SOCP, and SDP for the special-case ED selection problem that is commonly encountered in tree breeding.
We discussed the strength of the three conic relaxation problems,
and gave the theoretical bounds of the randomized algorithm that uses the SDP relaxation problem.
The fact that the theoretical bounds are not so sharp motivated us to implement
the steep-ascent method so that we can acquire a suitable solution for practical usage.
From the numerical results, we found that the steep-ascent method
with the SOCP relaxation was the most effective among the three conic
relaxation approaches, and that this outperformed the existing methods,
in particular, from the viewpoint of computation time.

One of further directions is to find a better theoretical bounds for the 
conic relaxation problems. In the discussions of this paper,
we mainly relied on the positive semidefiniteness of and non-negativity of 
the Wright numerator matrix $\A$. Since the specific values of the elements
in this matrix strongly relate to the pedigree of the candidate pool,
there is a possibility that we exploit such structures 
to tighten the theoretical bounds discussed in Lemma~\ref{le:SDP-R}.
However, we also need to reduce the computation time of the SDP relaxation problem
to make the SDP approach effective.

Another research direction is to minimize inbreeding depression~\cite{lindgren2009unequal}.
The objective function there is of form $(1- (\mbox{ID}) \x^T \A \x) \g^T \x$, where
ID is a constant that represents a regression slope.
Since the function is a cubic function with respect to the contribution $\x$,
this function is not even a convex function.
From the numerical results in this paper, 
we expect that a similar method to the steep-ascent method may perform well
for solving such a problem.
The minimization of the inbreeding depression will be an interesting problem
to researchers in the optimization field.

\setlength\baselineskip{1pt}%←Adjust here
\setlength{\itemsep}{0cm}
\setlength{\parskip}{0cm}
\bibliographystyle{abbrv}
\bibliography{reference} 

\end{document}